\documentclass[psamsfonts, reqno, 12pt]{amsart}
\usepackage{amsfonts,amsmath, amsthm, amssymb, latexsym, epsfig}
\usepackage{geometry}                
\usepackage{graphicx}
\usepackage{amssymb}
\usepackage{epstopdf}
\usepackage[all]{xy}
\usepackage{xspace}
\usepackage[mathscr]{eucal} 

	\advance \topmargin by -\headheight
	\advance \topmargin by -\headsep

	\evensidemargin \oddsidemargin
\RequirePackage{amsfonts}
\RequirePackage{amsmath}
\RequirePackage{amsthm}
\RequirePackage{amssymb}
\RequirePackage{latexsym}
\RequirePackage{epsfig}
\RequirePackage{verbatim}
\RequirePackage[mathscr]{eucal}
\RequirePackage[all]{xy}
	
	\newcommand{\sa}{\ensuremath{\operatorname{s.a.}}}
	\newcommand{\stex}[5]{\xymatrix{0\ar[r]&#1\ar[r]^{#2}&#3\ar[r]^{#4}&#5\ar[r]&0}}	
	\newcommand{\proj}{\ensuremath{\operatorname{\mathbb{P}}}}
	\newcommand{\setof}[2]{\ensuremath{\left\{ #1 \: : \: #2 \right\}}}

	\newcommand{\ftn}[3]{ #1 : #2 \rightarrow #3 }
	\newcommand{\limit}[3]{ \lim_{ #1 \to #2 } #3 }

	\newcommand{\norm}[1]{\left\|#1\right\|}

	\newcommand{\unit}[1]{1_{ #1 }}
	\newcommand{\unitize}[1]{\widetilde{#1}}

	\newcommand{\totalk}{\ensuremath{\underline{K} }\xspace} 
	
	\newcommand{\mc}[1]{\mathcal{#1}}
	
	\newcommand{\msf}[1]{\mathsf{#1}}
	\newcommand{\mrm}[1]{\mathrm{#1}}

	\newcommand{\mbb}[1]{\mathbb{#1}}
	
	\newcommand{\mbf}[1]{\mathbf{#1}}

	\newcommand{\Z}{\ensuremath{\mathbb{Z}}\xspace}
	\newcommand{\C}{\ensuremath{\mathbb{C}}\xspace}
	
	\newcommand{\R}{\ensuremath{\mathbb{R}}\xspace}
	\newcommand{\N}{\ensuremath{\mathbb{N}}\xspace}

	\newcommand{\kl}{\ensuremath{\mathit{KL}}\xspace}
	\newcommand{\kk}{\ensuremath{\mathit{KK}}\xspace}

	\newcommand{\Hom}{\ensuremath{\operatorname{Hom}}}

	\newcommand{\innerauto}{\ensuremath{\operatorname{Ad}}}
	\newcommand{\diag}{\ensuremath{\operatorname{diag}}}
	
	\newcommand{\id}{\ensuremath{\operatorname{id}}}

        \newcommand{\Extab}{\ensuremath{\operatorname{Ext}}_{\Z}^{1}}
	
	\newcommand{\pext}{\ensuremath{\operatorname{Pext}}_{\Z}^{1}}
	
	\newcommand{\ext}{\ensuremath{\operatorname{ext}}_{\Z}^{1}}

	\newcommand{\corona}[1]{\mathcal{Q}(#1)\xspace}

	\newcommand{\dirlim}{\displaystyle \lim_{\longrightarrow}}
	\newcommand{\invlim}{\displaystyle \lim_{\longleftarrow}}

	\newcommand{\cstar}{{$C \sp \ast$}\xspace}

	\theoremstyle{plain}
	\newtheorem{thm}{Theorem}[section]
	\newtheorem{lemma}[thm]{Lemma}
	\newtheorem{theorem}[thm]{Theorem}
	\newtheorem{proposition}[thm]{Proposition}
	\newtheorem{corollary}[thm]{Corollary}

	\theoremstyle{definition}
	\newtheorem{definition}[thm]{Definition}
	
	\newtheorem{remark}[thm]{Remark}
	
	\newtheorem{example}[thm]{Example}

	\newtheorem{question}[thm]{Question}
	
	\numberwithin{equation}{section}
	\numberwithin{figure}{section}

\begin{document}
	\title{The Automorphism group of a simple tracially AI algebra}
	\author{Ping Wong Ng}
	\address{University of Louisiana at Lafayette \\
	217 Maxim D. Doucet Hall \\
 	P.O.Box 41010 \\
	Lafayette, LA 70504-1010 USA }
	\email{png@louisiana.edu}
	\author{Efren Ruiz}
        \address{Department of Mathematics \\
        University of Toronto \\
        40 St. George St. \\
        Toronto, Ontario M5S 2E4 Canada}
        \email{eruiz@math.toronto.edu}
        \date{\today}
	


	\begin{abstract}
	The structure of the automorphism group of a simple TAI algebra is studied.  In particular, we show that $\frac{ \overline{ \mrm{Inn} }  ( A ) }{ \overline{ \mrm{Inn} }_{0} ( A ) }$ is isomorphic (as a topological group) to an inverse limit of discrete abelian groups for a unital, simple, AH algebra with bounded dimension growth.  Consequently, $\frac{ \overline{ \mrm{Inn} }  ( A ) }{ \overline{ \mrm{Inn} }_{0} ( A ) }$ is totally disconnected.  
	
Another consequence of our results is the following:  Suppose $A$ is the transformation group \cstar-algebra of a minimal Furstenberg transformation $( \mbb{T}^{n} , h_{n} )$ with a unique $h_{n}$-invariant probability measure on $\mbb{T}^{n}$.  Then the automorphism group of $A$ is an extension of a simple topological group by the discrete group $\mrm{ Aut } ( \totalk(A) )_{+,1}$.  
	\end{abstract}
        \maketitle

\section{Introduction}

In recent years, the program to classify amenable, separable \cstar-algebra, known as the Elliott Classification program, have proved to be quite useful in other fields of mathematics.  In particular, the Elliott Classification program has a played a role in finding interesting equivalence relations between two dynamical systems.  The amazing result of Giordano, Putnam, and Skau  \cite{GPS} involving minimal Cantor systems is such an example.  The Elliott Classification program helped to show that orbit equivalence (in a strong sense) between two minimal Cantor systems are characterized by their naturally associated pointed ordered groups.  More precisely, Giordano, Putnam, and Skau showed that two minimal Cantor systems are orbit equivalence (in a strong sense) if and only if their associated transformation group \cstar-algebras are isomorphic.  It was shown that for every minimal Cantor system, the associated transformation group \cstar-algebra is a simple AT algebra with real rank zero.  Then, by Elliott's classification of AT algebras with real rank zero \cite{at}, these transformation group \cstar-algebras are completely classified by their ordered $K$-groups together with the class of the unit, which are the naturally associated pointed ordered group of a minimal Cantor system.   

One of the goal of this paper is to show that the Elliott Classification program can be used to study the structure of the automorphism group of a \cstar-algebra $A$, where $A$ is in a class of \cstar-algebras $\mc{C}$ in which every \cstar-algebra in $\mc{C}$ is classified via $K$-theory.  In particular, we will study the structure of the automorphism group of a \cstar-algebra $A$ arising from a dynamical system.  This class of \cstar-algebras includes the irrational rotation algebras, the \cstar-algebras considered by Giordano, Putnam, and Skau, the simple Cuntz-Krieger algebras, and the simple noncommutative tori.   

Set $\overline{ \mrm{Inn} } ( A )$ to be the group of approximately inner automorphisms of $A$ (equip with the topology of pointwise convergence) and $\overline{ \mrm{ Inn } }_{0} ( A )$ to be closure of the group of inner automorphisms of $A$ whose implementing unitaries is connected to $1_{A}$ via a continuous path of unitaries.  For any groups $G$ and $H$ and for any $u$ in $G$, denote the subgroup of $H$ consisting of all elements $h$ of $H$ such that there exists a group homomorphism from $G$ to $H$ which sends $u$ to $h$ by $H_{u} ( G, H )$.

\bigskip

\noindent \textbf{Theorem \ref{thm:inniso}.} \ \emph{Let $A$ be a unital, separable \cstar-algebra satisfying Property (C).  Suppose the natural map from $U(A) / U_{0} (A)$ to $K_{1} (A)$ is an isomorphism.  Then for every increasing sequence $\{ G_{n} \}_{n = 1}^{ \infty }$ of finitely generated subgroups of $K_{0} ( A )$ containing $[ 1_{A} ]$, 
\begin{equation*}
\frac{ \overline{ \mrm{Inn} } ( A ) }{ \overline{ \mrm{Inn} }_{0} ( A ) } \cong \invlim \frac{ K_{1} ( A ) }{ H_{ [ 1_{A} ] } ( G_{n} , K_{1} ( A ) }
\end{equation*}
}
as topological groups.

\bigskip

\noindent \textbf{Theorem \ref{thm:propC}.} \ \emph{Let $A$ be a unital, separable, amenable, simple \cstar-algebra satisfying the Universal Coefficient Theorem of Rosenberg and Schochet \cite{uct}.  If $A$ is either a purely infinite \cstar-algebra or a tracially AI algebra, then $A$ satisfies Property (C). }

\bigskip
 
Consequently, $\frac{ \overline{ \mrm{Inn } } ( A ) }{ \overline{ \mrm{Inn} }_{0} ( A ) }$ is a totally disconnected topological group.  We also give criterions for when $\frac{ \overline{ \mrm{Inn } } ( A ) }{ \overline{ \mrm{Inn} }_{0} ( A ) }$ is a discrete group, compact group, or locally group.

Our results generalize the results of Elliott and R{\o}rdam in \cite{EllRorAuto}.  In their paper, Elliott and R{\o}rdam showed that $\frac{ \overline{ \mrm{Inn} } ( A ) }{ \overline{ \mrm{Inn} }_{0} ( A ) }$ is a totally disconnected group whenever $A$ is a unital, simple, AT algebra with real rank zero.  Elliott and R{\o}rdam showed that $\frac{ \overline{ \mrm{Inn} } ( A ) }{ \overline{ \mrm{Inn } }_{0} ( A ) }$ is isomorphic (as a topological group) to the inverse limit of the discrete groups $K_{1} ( A ) / n K_{1} ( A )$, where $n$ belongs to a directed set of positive integers that divide $ [ 1_{A} ] $ in $K_{0} ( A )$.  It turns out that there are examples in which the inverse limit of discrete groups introduced by Elliott and R{\o}rdam is not isomorphic to the inverse limit of discrete groups introduced in this paper.  Take for example, $K_{0} ( A ) = \Z / n \Z$ and $K_{1} ( A ) = \Z$.  Then, the inverse limit of Elliott and R{\o}rdam is a torsion group and the inverse limit considered in this paper is isomorphic to $\Z$, given the discrete topology.    

We apply our results to the transformation group \cstar-algebra associated to a minimal Furstenberg transformation with a unique invariant probability measure.  More precisely, define $h_{n}$ from $\mbb{T}^{n}$ to $\mbb{T}^{n}$ to be the inverse of the homeomorphism
\begin{equation*}
( \zeta_{1} , \zeta_{2} , \dots, \zeta_{n} ) \mapsto ( e^{ 2 \pi i \theta } \zeta_{1} , f_{1} ( \zeta_{1} ) \zeta_{2} , f_{2} ( \zeta_{1 } , \zeta_{2} ) \zeta_{3} , \dots , f_{n-1} ( \zeta_{1} , \dots , \zeta_{n-1}) \zeta_{n} ).
\end{equation*}    
where each $f_{i}$ is a continuous function from $\mbb{T}^{i}$ to $\mbb{T}$.  If each $f_{i}$ is Lipschitz, then Furstenberg in \cite{furstenberg} showed that $( \mbb{T}^{n} , h_{n} )$ is minimal and has a unique $h_{n}$-invariant probability measure.  Let $A$ denote the transformation group \cstar-algebra associated to $( \mbb{T}^{n} , h_{n} )$.

\bigskip

\noindent \textbf{Theorem \ref{thm:furstenberg}.} \ \emph{If $A$ is as above, then $A$ is a separable, unital, amenable, simple, tracially AF algebra which satisfies the Universal Coefficient Theorem of Rosenberg and Schochet \cite{uct}.  Moreover, $\mrm{Aut} ( A )$ fits into the following short exact sequence:
\begin{equation*}
\{ 1 \} \to \overline{ \mrm{Inn } } ( A ) \to \mrm{Aut} ( A ) \to \mrm{Aut} ( \totalk(A) )_{+,1} \to \{ 0 \}.
\end{equation*}  
and $\overline{ \mrm{Inn} } ( A ) = \overline{ \mrm{Inn} }_{0} ( A )$.   Consequently, $\overline{ \mrm{Inn} } ( A )$ is a simple topological group.
}
 
\bigskip
   
The paper is organized as follows.  In Section \ref{innerauto}, we define Property (C) and show that if $A$ is a unital, separable \cstar-algebra satisfying Property (C), then $\frac{ \overline{ \mrm{Inn} } ( A ) }{ \overline{ \mrm{Inn} }_{0} ( A ) }$ is a totally disconnected group (see Theorem \ref{thm:inn}).  If, in addition, the canonical map from $U(A)/U_{0}(A)$ to $K_{1}(A)$ is an isomorphism, then $\frac{ \overline{ \mrm{Inn} } ( A ) }{ \overline{ \mrm{Inn} }_{0} ( A ) }$ is isomorphic (as a topological group) to an inverse limit of discrete abelian groups.  Section \ref{cstarPropC} is devoted to showing that $A$ satisfies Property (C) whenever $A$ is a unital, separable, amenable, simple \cstar-algebra satisfying the Universal Coefficient Theorem, and $A$ is either a purely infinite \cstar-algebra or a tracially AI algebra (Theorem \ref{thm:propC}).  In the last section, we apply our results to study the automorphism group of $A$, where $A$ is the transformation group \cstar-algebra associated to a minimal Furstenberg transformation $h$ with a unique $h$-invariant measure on $\mbb{T}^{n}$.  We also give an explicit description of $\mrm{Aut} ( \totalk ( A ) )_{+,1}$ when $n$ is equal to $2$.

\section{The closure of the group of inner automorphisms}\label{innerauto}

\noindent We will use the following conventions throughout the paper.

\begin{enumerate}
\item Let $A$ be a separable \cstar-algebra. $A$ is said to satisfy the \emph{Universal Coefficient Theorem} if for any $\sigma$-unital \cstar-algebra $B$, the sequence
\begin{equation*}
0 \to \Extab(K_{*}(A), K_{* + 1} (B) ) \overset{ \delta }{\to}   \kk^{0}(A,B) \overset{ \Gamma }{ \to } \Hom(K_{*}(A), K_{*}(B)) \to 0
\end{equation*}
is exact. 

\item  Denote the set of all equivalence classes of locally trivial extensions in $\Extab( G, H )$ by $\pext( G, H )$.  It is easy to see that $\pext( G, H )$ is a subgroup of $\Extab( G, H )$.  Denote the quotient by $\ext( G, H )$.  Let $A$ be a \cstar-algebra that satisfies the Universal Coefficient Theorem and let $B$ be a $\sigma$-unital \cstar-algebra.  Then set $\kl^{0}(A,B) = \kk(A,B)/ \pext(K_{*}(A), K_{*+1}(B))$ and $\kl^{1}(A,B) = \kl^{0}(A, SB)$ (see R\o rdam \cite{defKL}).  We will also denote $\kl^{0}(A,B)$ by $\kl ( A, B )$.  Since $A$ satisfies the Universal Coefficient Theorem, the following are exact sequences:
\begin{equation*}
0 \to \ext(K_{*}(A), K_{* + 1} (B) ) \overset{ \delta }{ \to }  \kl^{0}(A,B) \overset{ \Gamma }{ \to } \Hom(K_{*}(A), K_{*}(B)) \to 0
\end{equation*}
\begin{equation*}
0 \to \ext(K_{*}(A), K_{* } (B) ) \overset{ \delta }{ \to }  \kl^{1}(A,B) \overset{ \Gamma }{ \to } \Hom(K_{*}(A), K_{*+1}(B)) \to 0.
\end{equation*}

\item Let $A$ be a unital \cstar-algebra.  Denote the tracial state space of $A$ by $T(A)$.  The space of all affine continuous map from $T(A)$ to $\R$ will be denoted by $\mrm{Aff} ( T(A) )$.  For every projection $p$ in $\msf{M}_{n} (A)$, define $f_{p}$ in $\mrm{Aff} ( T(A) )$ to be the function which sends $\tau$ in $T(A)$ to $\tau ( p )$.  The map from $K_{0} ( A )$ to $\mrm{Aff} ( T ( A ) )$ which sends $[ p ] $ to the map $f_{ p }$ will be denoted by $\rho_{ A }$. 

\item  Let $G$ and $H$ be groups and let $g$ be an element of $G$.  Then $H_{g} ( G, H )$ will denote the subgroup of $H$ consisting of all elements $h$ of $H$ such that there exists a homomorphism $\alpha$ from $G$ to $H$ with $\alpha ( g ) = h$.  If $G = K_{0} ( A )$, $H = K_{0} ( B )$, $g = [ 1_{A} ]$, and $\corona{B}$ satisfies a certain infinite property, then Lin showed in \cite{unitaryequiv} that two unital essential extension $\tau_{1}$ and $\tau_{2}$ from $A$ to $\corona{B}$ that are unitarily equivalent are strongly unitarily equivalent if and only if the implement unitary can be chosen to be in $H_{ [ 1_{A} ] } ( K_{0} ( A ) , K_{0} ( B ) )$.
\end{enumerate}

We will now explain the notation used in Lemma \ref{lem:invlim} below and in the process prove Lemma \ref{lem:invlim}.

Let $A$ be a separable unital \cstar-algebra and let $\{ G_{n} \}_{ n = 1}^{ \infty }$ be an increasing sequence of finitely generated subgroups of $[ 1_{A} ]$ with each $G_{n}$ containing $[ 1_{A} ]$.  By the choice of $\{ G_{n} \}_{ n = 1}^{ \infty }$, there exists a sequence of surjections
\begin{equation*}
\frac{ K_{1} ( A ) }{ H_{ [ 1_{A} ] } ( G_{1}, K_{1} ( A ) ) } \leftarrow \frac{ K_{1} ( A ) }{ H_{ [ 1_{A} ] } ( G_{2}, K_{1} ( A ) ) } \leftarrow \cdots,
\end{equation*}  
which gives an inverse limit group
\begin{equation*}
\invlim \frac{ K_{1} ( A ) }{ H_{ [ 1_{A} ] } ( G_{n}, K_{1} ( A ) ) }.
\end{equation*}

Equip $\frac{ K_{1} ( A ) }{ H_{ [ 1_{A} ] } ( G_{n} , K_{1} ( A ) ) }$ with the discrete topology and give the inverse limit group the natural inverse limit topology.  Denote the natural projection 
\begin{equation*}
\invlim \frac{ K_{1} ( A ) }{ H_{ [ 1_{A} ] } ( G_{n}, K_{1} ( A ) ) } \to \frac{ K_{1} ( A ) }{ H_{ [ 1_{A} ] } ( G_{n}, K_{1} ( A ) )  }
\end{equation*}
by $\pi_{n}$.  Then a sequence $\{ x_{k} \}_{ k = 1}^{ \infty }$ converges to $x$ in the inverse limit if and only if the sequence $\{ \pi_{n} ( x_{k} ) \}_{ k = 1}^{ \infty }$ converges to $\pi_{n} ( x )$ in $\frac{ K_{1} ( A ) }{ H_{ [ 1_{A} ] } ( G_{n}, K_{1} ( A ) )  }$ for all $n$ in $\N$ if and only if for all $n$ in $\N$, there exists $k(n)$ in $\N$ such that $\pi_{n} ( x_{k} ) = \pi_{n} ( x )$ for all $k \geq k(n)$.

Since the quotient maps from $K_{1} ( A )$ to $\frac{ K_{1} ( A ) }{ H_{ [ 1_{A} ] } ( G_{n}, K_{1} ( A ) ) }$ commute with the surjections in the inverse limit, they lift to a map
\begin{equation*}
K_{1} ( A ) \to \invlim \frac{ K_{1} ( A ) }{ H_{ [ 1_{A} ] } ( G_{n}, K_{1} ( A ) ) }.
\end{equation*}
For each $g$ in $K_{1} ( A )$ denote the image of $g$ in $\invlim \frac{ K_{1} ( A ) }{ H_{ [ 1_{A} ] } ( G_{n}, K_{1} ( A ) ) }$ by $\check{g}$.  

From the above paragraphs, we have the follow result. 

\begin{lemma}\label{lem:invlim}
Let $A$ and $\{ G_{n} \}_{ n = 1}^{ \infty }$ be as in the above paragraphs.  Then the image of $K_{1} ( A )$ in $\invlim \frac{ K_{1} ( A ) }{ H_{ [ 1_{A} ] } ( G_{n}, A ) }$ is dense.  If $\{ g_{k} \}_{ k = 1}^{ \infty }$ is a sequence in $K_{1}( A )$, then the sequence $\{ \check{g}_{k} \}_{ k = 1}^{ \infty }$ converges to zero in $\invlim \frac{ K_{1} ( A ) }{ H_{ [ 1_{A} ] } ( G_{n}, A ) }$ if and only if for each $n$ in $\N$, there exists $k(n)$ in $\N$ such that for all $k \geq k(n)$, we have that $g_{k}$ is an element of $H_{ [ 1_{A} ] } ( G_{n} , K_{1} ( A ) )$.  
\end{lemma} 

Recall that $\overline{ \mrm{Inn} } ( A )$ is the closure of $\setof{ \innerauto( u ) }{ u \in U(A) }$ under the topology induced by pointwise convergence and $\overline{ \mrm{Inn} }_{0} ( A )$ is the normal subgroup of $\overline{ \mrm{Inn} } ( A )$ which is the closure of $\setof{ \innerauto ( u ) }{ u \in U_{0} ( A ) }$. Thus, there exists a homomorphism from $U ( A ) / U_{0} ( A )$ to $\frac{ \overline{ \mrm{Inn} } ( A ) }{ \overline{ \mrm{ Inn } }_{0} ( A ) }$ such that the diagram
\begin{equation}\label{diagramInn}
\vcenter{\xymatrix{
\{1\} \ar[r] & U_{0} ( A ) \ar[r] \ar[d]_{ \mrm{Ad} } & U(A) \ar[r] \ar[d]_{ \mrm{Ad} } & U( A ) / U_{0} ( A ) \ar[r] \ar[d] & \{ 1 \} \\
\{ 1 \} \ar[r] & \overline{ \mrm{Inn} }_{0} ( A ) \ar[r] & \overline{ \mrm{Inn} } ( A ) \ar[r] & \frac{ \overline{ \mrm{Inn} } ( A ) }{ \overline{ \mrm{Inn} }_{0} ( A ) } \ar[r] & \{ 1 \}
}}
\end{equation} 
is commutative.  If $\frac{ \overline{ \mrm{Inn} } ( A ) }{ \overline{ \mrm{Inn} }_{0} ( A ) }$ is given the quotient topology, then $\frac{ \overline{ \mrm{Inn} } ( A ) }{ \overline{ \mrm{Inn} }_{0} ( A ) }$ is a complete topological group.  

For every $g$ in $U( A ) / U_{0} ( A )$, denote the image of $g$ in $\frac{ \overline{ \mrm{Inn} } ( A ) }{ \overline{ \mrm{Inn} }_{0} ( A ) }$ via the map from $U(A)/U_{0}(A)$ to $\frac{ \overline{ \mrm{Inn} } ( A ) }{ \overline{ \mrm{Inn} }_{0} ( A ) }$ in the above diagram by $\hat{g}$.    

We will show that there exists a continuous homomorphism from $\frac{ \overline{ \mrm{Inn} } ( A ) }{ \overline{ \mrm{Inn} }_{0} ( A ) }$ to $\invlim \frac{ K_{1} ( A ) }{ H_{ [1_{A}] } ( G_{n} , K_{1} ( A ) ) }$ for any increasing sequence $\{ G_{n} \}_{ n = 1}^{ \infty }$ of finitely generated subgroup of $K_{0} ( A )$ with each $G_{n}$ containing $[ 1_{A} ]$ and whose union is equal to $K_{0} ( A )$.  In order to prove this, we first need the following result.  The proof of the lemma which is similar to Lemma 4.2 of \cite{EllRorAuto}.  In fact, the proof carries over to our case.  Since our statement is slightly different we prove it here.   

\begin{lemma}\label{lem:stronginn}
Let $A$ be a unital, separable \cstar-algebra.  The image of $U(A) / U_{0} (A)$ (via the map $g \mapsto \hat{g}$) is dense in $\frac{ \overline{ \mrm{Inn} } ( A ) }{ \overline{ \mrm{Inn} }_{0} ( A ) }$.  Moreover, if $\{ g_{k} \}_{ k = 1}^{\infty}$ is a sequence in $U( A ) / U_{0} (A)$, then the sequence $\{ \widehat{g_{k}} \}_{ k = 1}^{ \infty }$ converges to zero in $\frac{ \overline{ \mrm{Inn} } ( A ) }{ \overline{ \mrm{Inn} }_{0} ( A ) }$ if and only if there exists a sequence of unitaries $\{ u_{k} \}_{ k = 1}^{ \infty }$ in $A$ such that
\begin{enumerate}
\item for all $k$, $g_{k} = [ u_{k} ]$ in $U(A) / U_{0} ( A )$ and 
\item $\limit{k}{\infty}{ \norm{ u_{k} a - a u_{k} } } = 0$ for all $a$ in $A$.
\end{enumerate}
\end{lemma}

\begin{proof}
The density of the image of $U(A) / U_{0} (A)$ follows from Diagram (\ref{diagramInn}).  Now note that the sequence $\{ \widehat{g_{k}} \}_{ k = 1}^{ \infty}$ converges to zero in $\frac{ \overline{ \mrm{Inn} } ( A ) }{ \overline{ \mrm{Inn} }_{0} ( A ) }$ if and only if the sequence lifts to a sequence in $\overline{ \mrm{Inn} } ( A )$ that converges to $\id_{A}$ in $\overline{ \mrm{Inn} } ( A )$.  

Suppose the sequence $\{ \widehat{g_{k}} \}_{ k = 1}^{ \infty}$ converges to zero in $\frac{ \overline{ \mrm{Inn} } ( A ) }{ \overline{ \mrm{Inn} }_{0} ( A ) }$.  Let $\alpha_{k}$ be a lifting of $\widehat{g_{k}}$ such that $\{ \alpha_{k} \}_{ k = 1}^{ \infty }$ converges to $\id_{A}$ in $\overline{ \mrm{Inn} } ( A )$ and let $w_{k}$ be an element of $U(A)$ which is a lifting of $g_{k}$.  Since $\innerauto ( w_{k} )$ is a lifting of $\widehat{ g_{k} }$ and $\alpha_{k}$ is a lifting of $\widehat{ g_{k} }$, we have that $\alpha_{k} = \beta_{k} \innerauto( w_{k} )$ and $\beta_{k}$ is an element of $\overline{ \mrm{Inn} }_{0} ( A )$.  Since the image of $U_{0} ( A )$ via the map $\mrm{Ad}$ is dense in $\overline{ \mrm{Inn} }_{0} ( A )$ and the topological group $\overline{ \mrm{Inn} }_{0} ( A )$ is first countable, there exists a sequence of unitaries $\{ v_{k} \}_{ k = 1}^{ \infty }$ in $U_{0} (A)$ such that the sequence $\{ \beta_{k} \gamma_{k}^{-1} \}_{ k = 1}^{ \infty }$ converges to $\id_{A}$ in $\overline{ \mrm{Inn} }_{0} ( A )$ where $\gamma_{k} = \innerauto( v_{k} )$.  Hence, the sequence $\{ \gamma_{k} \innerauto( w_{k} ) \}_{ k = 1}^{ \infty }$ converges to $\id_{A}$ in $\overline{ \mrm{Inn} } (A)$ and $\gamma_{k} \innerauto ( w_{k} ) = \innerauto ( v_{k} w_{k} )$ is a lifting of $\widehat{g_{k}}$.  

Set $u_{k} = v_{k} w_{k}$.  Then 
\begin{enumerate}
\item $[ u_{k} ] = g_{k}$ in $U(A) / U_{0} (A)$ for all $k$ and 
\item $\limit{k}{ \infty }{ \norm{ u_{k} a - a u_{k} } } = 0$ for all $a$ in $A$. 
\end{enumerate}

Suppose there exists a sequence of unitaries $\{ u_{k} \}_{ k = 1}^{ \infty }$ in $A$ such that
\begin{enumerate}
\item $[ u_{k} ] = g_{k}$ in $U(A) / U_{0} (A)$ for all $k$ and 
\item $\limit{k}{ \infty }{ \norm{ u_{k} a - a u_{k} } } = 0$ for all $a$ in $A$. 
\end{enumerate}
Then it easy to see that $\innerauto ( u_{k} )$ is a lifting of $\widehat{g_{k}}$ and the sequence $\{ \innerauto( u_{k} ) \}_{ k = 1}^{ \infty }$ converges to $\id_{A}$ in $\overline{ \mrm{Inn} } ( A )$.
\end{proof}

We will also need the following fact from group theory.  Since we have not been able to find a reference to the result below, for the convenience of the reader, we prove it here.

\begin{lemma}\label{lem:locallytrivial}
Let $\{ G_{n} \}_{ n = 1}^{ \infty }$ be a sequence of abelian groups and $G$ be an abelian group.  Suppose $\alpha$ from $G$ to $\frac{ \prod_{ n =1}^{ \infty } G_{n} }{ \bigoplus_{ n = 1}^{ \infty} G_{n} }$ is a group homomorphism.  Denote the quotient homomorphism from $\prod G_{n}$ to $\frac{ \prod G_{n} }{ \bigoplus_{ n = 1}^{ \infty} G_{n} }$ by $\epsilon$.  Then for every finitely generated subgroup $H$ of $G$, there exists a group homomorphism $\beta_{H}$ from $H$ to $\prod_{ n = 1}^{ \infty } G_{n}$ such that $\alpha \vert_{ H } = \epsilon \circ \beta_{H}$.  

Consequently, the short exact sequence 
\begin{equation*}
\xymatrix{
0 \ar[r] & \bigoplus_{ n = 1}^{ \infty } G_{n} \ar[r] & \prod_{ n = 1}^{ \infty } G_{n} \ar[r] & \frac{ \prod_{ n =1}^{ \infty } G_{n} }{ \bigoplus_{ n = 1}^{ \infty} G_{n} } \ar[r] & 0
}
\end{equation*}
is locally trivial, i.e. restricting the above short exact sequence to any finitely generated subgroup of $\frac{ \prod_{ n =1}^{ \infty } G_{n} }{ \bigoplus_{ n = 1}^{ \infty} G_{n} }$ the short exact sequence obtained by this process is a split exact sequence.
\end{lemma}

\begin{proof}
Note that it is enough to show that for every $x$ in $H$ with finite order $k$, there exists $y$ in $\prod_{ n = 1}^{ \infty } G_{n}$ such that $k y = 0$ and $\alpha ( x ) = \epsilon ( y )$.  Let $x$ be an element of $H$ with finite order $k$.  Note that $k \alpha ( x ) = 0$.  Take any lifting $a = \{ a_{n} \}_{ n = 1}^{ \infty }$ of $\alpha ( x )$.  Since $k \alpha ( x ) = 0$, we have that $ k a $ is an element of $\bigoplus_{ n = 1}^{ \infty } G_{n}$.  Hence, there exists $n_{0}$ such that for all $n$ greater than $n_{0}$, we have that $k a_{n} = 0$.  Define $y = \{ y_{n} \}_{ n = 1}^{ \infty }$ in $\prod_{ n = 1}^{ \infty } G_{n}$ by $y_{n} = 0$ if $n$ is less than or equal to $n_{0}$ and $y_{n} = a_{n}$ if $n$ is greater than $n_{0}$.  Then $k y = 0$ and $\epsilon ( y ) = \alpha ( x )$. 
\end{proof}
    
Let $A$ be a unital \cstar-algebra.  For notational convenience, for every element $g$ of $U(A) / U_{0} (A)$ we will again denote the image of $g$ under the canonical map from $U(A) / U_{0} (A)$ to $K_{1} (A)$ by $g$.    
    
\begin{theorem}\label{thm:surjinn}
Let $A$ be a separable, unital \cstar-algebra.  Let $\{ G_{n} \}_{ n = 1}^{ \infty }$ be an increasing sequence of finitely generated subgroups of $K_{0} ( A )$ whose union is $K_{0} ( A )$ and for each $n$, $G_{n}$ contains $[ 1_{A} ]$.  Then the map $\mu$ which sends $\hat{g}$ to $\check{g}$ extends uniquely to a continuous group homomorphism from $\frac{ \overline{ \mrm{Inn} } ( A ) }{ \overline{ \mrm{Inn} }_{0} ( A ) }$ to $\invlim \frac{ K_{1} ( A ) }{ H_{ [1_{A}] } ( G_{n} , K_{1} ( A ) ) }$.  We will denote this extension by $\mu$.
\end{theorem}

\begin{proof}
Note that the topological groups $\frac{ \overline{ \mrm{Inn} } ( A ) }{ \overline{ \mrm{Inn} }_{0} ( A ) }$ and  $\invlim \frac{ K_{1} ( A ) }{ H_{ [ 1_{A} ] } ( G_{n} , K_{1} ( A ) ) }$ are complete.  Therefore, by Lemma \ref{lem:stronginn} $\frac{ \overline{ \mrm{Inn} } ( A ) }{ \overline{ \mrm{Inn} }_{0} ( A ) }$ is the completion of the image of $U( A ) / U_{0} (A)$ and by Lemma \ref{lem:invlim}, $\invlim \frac{ K_{1} ( A ) }{ H_{ [ 1_{A} ] } ( G_{n} , K_{1} ( A ) ) }$ is the completion of the image of $K_{1} ( A )$.  Hence, to prove the theorem it is enough to show that for every sequence $\{ g_{k} \}_{ k = 1}^{ \infty }$ in $U ( A ) / U_{0} (A)$, the sequence $\{ \hat{ g_{k} } \}_{ k = 1}^{ \infty }$ converges to zero in $\frac{ \overline{ \mrm{Inn} } ( A ) }{ \overline{ \mrm{Inn} }_{0} ( A ) }$ implies that the sequence $\{ \check{g}_{k} \}_{ k = 1}^{ \infty }$ converges to zero in $\invlim \frac{ K_{1} ( A ) }{ H_{ [ 1_{A} ] } ( G_{n} , K_{1} ( A ) ) }$.  

Suppose $\{ \hat{g}_{k} \}_{ k = 1 }^{ \infty }$ converges to zero in $\frac{ \overline{ \mrm{Inn} } ( A ) }{ \overline{ \mrm{Inn} }_{0} ( A ) }$.  By Lemma \ref{lem:stronginn}, there exists a sequence of unitaries $\{ u_{k} \}_{ k = 1}^{ \infty }$ in $A$ such that 
\begin{enumerate}
\item for all $k$ in $\N$, we have that $[ u_{k} ] = g_{k}$ in $U( A ) / U_{0} (A)$ and 
\item $\limit{ k }{ \infty }{ \norm{ u_{k} a - a u_{k} } } = 0$ for all $a$ in $A$.
\end{enumerate} 

Denote the \cstar-direct product of infinite copies of $A$ by $\ell^{ \infty } ( A )$ and denote the ideal of $\ell^{ \infty } ( A )$ consisting of all sequences $\{ x_{k} \}_{ k = 1}^{ \infty }$ in $\ell^{ \infty } ( A )$ which approaches to zero by $c_{0} ( A )$.  Denote the quotient of $\ell^{ \infty } ( A )$ with $c_{0} ( A )$ by $q_{ \infty } ( A )$ and denote the quotient map from $\ell^{ \infty } ( A )$ to $q_{ \infty } ( A )$ by $\pi_{ \infty }$.  Define $\varphi$ from $A \otimes C(S^{1})$ to $q_{ \infty } ( A )$ by $\varphi ( a \otimes 1_{ C(S^{1}) } ) = \pi_{ \infty } ( \{ a \} )$ and $\varphi ( 1_{A} \otimes z ) = \pi_{ \infty } ( \{ u_{k} \} )$.  Since $\limit{ k }{ \infty }{ \norm{ u_{k} a - a u_{k} } } = 0$ for all $a$, we have that $\varphi$ is a unital $*$-homomorphism.

Note that in the diagram
\begin{equation*}
\xymatrix{ 0 \ar[r] & K_{1} ( c_{0} ( A ) ) \ar[r] \ar[d] & K_{1} ( \ell^{\infty} (A) ) \ar[r] \ar[d] & K_{1} ( q_{ \infty} ( A ) ) \ar[r] \ar[d] & 0 \\
		0 \ar[r] & \bigoplus_{ k = 1}^{ \infty } K_{1} ( A ) \ar[r] & \prod_{ k = 1}^{ \infty } K_{1} ( A ) \ar[r] & \frac{ \prod_{k = 1}^{\infty} K_{1} ( A ) }{ \bigoplus_{k = 1}^{ \infty } K_{1} ( A ) } \ar[r] & 0
		}
\end{equation*}
the rows are exact and the middle vertical map is define by $[ \{ w_{k} \} ] \mapsto \{ [ w_{k} ] \}$.  Hence, by composing $K_{1} ( \varphi )$ with the vertical map in the third column of the above diagram, we get a homomorphism $\alpha$ from $K_{1} ( A \otimes C(S^{1}) )$ to $\frac{ \prod_{k = 1}^{\infty} K_{1} ( A ) }{ \bigoplus_{k = 1}^{ \infty } K_{1} ( A ) }$.  

Note that there exists a homomorphism $\beta$ from $K_{0} ( A )$ to $K_{1} ( A \otimes C(S^{1}) )$ such that $\beta$ maps $[ 1_{A} ] $ to $[ 1_{A} \otimes z ]$.  If $G$ is any finitely generated subgroup of $K_{0} ( A )$ containing $[ 1_{A} ]$, by Lemma \ref{lem:locallytrivial} $\alpha \circ \beta \vert_{G}$ lifts to a homomorphism from $G$ to $\prod_{ k = 1}^{ \infty } K_{1} (A)$.  Note that $\beta \circ \alpha$ maps $[ 1_{A} ]$ to the image of $\{ [ u_{k} ] \}_{ k =1}^{ \infty }$ in $\frac{ \prod_{k = 1}^{\infty} K_{1} ( A ) }{ \bigoplus_{k = 1}^{ \infty } K_{1} ( A ) }$.  

For each $n$, let $\gamma = \{ \gamma_{k} \}_{ k = 1}^{ \infty }$ from $G_{n}$ to $\prod_{ k = 1}^{ \infty } K_{1} ( A )$ be a lifting of $\alpha \circ \beta \vert_{ G_{n} }$.  Then $\{ \gamma_{k} ( [ 1_{A} ] ) - [ u_{k} ] \}_{ k = 1}^{ \infty }$ is an element of $\oplus_{ k = 1}^{ \infty } K_{1} ( A )$.  Hence, there exists $k(n)$ such that for all $k \geq k(n)$, we have that $\gamma_{k} ( [ 1_{A} ] ) = [ u_{k} ] = g_{k}$ in $K_{1} ( A )$.  This implies that for all $k \geq k(n)$, we have that $g_{k}$ is an element of $H_{ [ 1_{A} ] } ( G_{n} , K_{1} ( A ) )$.  By Lemma \ref{lem:invlim}, the sequence $\{ \check{g}_{k} \}_{ k = 1}^{ \infty }$ converges to zero in $\invlim \frac{ K_{1} ( A ) }{ H_{ [ 1_{A} ] } ( G_{n} , K_{1} ( A ) ) }$.  The uniqueness of $\mu$ is clear.  
\end{proof}
 
If $A$ is a unital separable \cstar-algebra satisfying Property (C) (defined below), then we will show that $\mu$ (see Theorem \ref{thm:surjinn}) is an injective map between topological groups.

\begin{definition}\label{def:propertyC}
A unital \cstar-algebra $A$ is said to satisfy Property (C) if for every $\epsilon$ in $\R_{ > 0}$ and finite subset $\mc{F}$ of $A$, there exists a finitely generated subgroup $G_{0}$ of $K_{0} ( A )$ containing $[ 1_{A} ]$ such that the following holds:  for every unitary $u$ in $A$ with $[u]$ an element of $H_{ [ 1_{A} ] } ( G_{0}, K_{1} ( A ) )$, there exists a unitary $w$ in $A$ such that
\begin{enumerate}
\item $\norm{ w a - a w } < \epsilon$ for all $a$ in $ \mc{F}$ and
\item $[ w ] = [ u ]$ in $U ( A ) / U_{0} ( A )$.
\end{enumerate}

Note that for every finitely generated subgroup $G$ of $K_{0} ( A )$ containing $G_{0}$, we have that $H_{ [ 1_{A} ] } ( G , K_{1} ( A ) )$ is a subgroup of $H_{ [ 1_{A} ] } ( G_{0} , K_{1} ( A ) )$.  Therefore, if $u$ is a unitary in $A$ with $[ u ]$ an element of $H_{ [ 1_{A} ] } ( G , K_{1} ( A ) )$, then there exists a unitary $w$ in $A$ such that 
\begin{enumerate}
\item $\norm{ w a - a w } < \epsilon$ for all $a$ in $\mc{F}$ and 
\item $[ w ] = [ u ]$ in $U ( A ) / U_{0} (A) $.
\end{enumerate}
\end{definition}

\begin{theorem}\label{thm:inn}
Let $A$ be a separable, unital \cstar-algebra that satisfies Property (C).  Then for any increasing sequence $\{ G_{n} \}_{n = 1}^{ \infty }$ of finitely generated subgroups of $K_{0} ( A )$ whose union is $K_{0} ( A )$ and each $G_{n}$ contains $[ 1_{A} ]$, the continuous map $\mu$ defined in Theorem \ref{thm:surjinn} from $\frac{ \overline{ \mrm{Inn} } ( A ) }{ \overline{ \mrm{Inn} }_{0} ( A ) }$ to $\invlim \frac{ K_{1} ( A ) }{ H_{ [ 1_{A} ] } ( G_{n} , K_{1} ( A ) ) }$ is injective. 

Consequently, $\frac{ \overline{ \mrm{Inn} } ( A ) }{ \overline{ \mrm{Inn} }_{0} ( A ) }$ is totally disconnected.
\end{theorem}

\begin{proof}
Note that to prove the theorem it is enough to show that for every sequence $\{ g_{k} \}_{ k = 1}^{ \infty }$ in $U ( A ) / U_{0} (A)$, the sequence $\{ \widehat{ g_{k} } \}_{ k = 1}^{ \infty }$ converges to zero in $\frac{ \overline{ \mrm{Inn} } ( A ) }{ \overline{ \mrm{Inn} }_{0} ( A ) }$ if and only if the sequence $\{ \check{g}_{k} \}_{ k = 1}^{ \infty }$ converges to zero in $\invlim \frac{ K_{1} ( A ) }{ H_{ [ 1_{A} ] } ( G_{n} , K_{1} ( A ) ) }$.  By Theorem \ref{thm:surjinn}, the only if direction holds since $\mu$ is continuous.  

Suppose that $\{ \check{g}_{k} \}_{ k = 1}^{ \infty }$ converges to zero in $\invlim \frac{ K_{1} ( A ) }{ H_{ [ 1_{A} ] } ( G_{n} , K_{1} ( A ) ) }$.  By Lemma \ref{lem:invlim}, for each $n$ in $\N$, there exists $k(n)$ in $\N$ such that for all $k \geq k(n)$ we have that $g_{k}$ is an element of $H_{ [ 1_{A} ] } ( G_{n} , K_{1} ( A ) )$.  Let $\{ \mc{F}_{m} \}_{ m = 1}^{ \infty }$ be an increasing sequence of finite subsets of $A$ whose union is dense in $A$.   We will show that there exist a sequence of unitaries $\{ u_{k} \}_{ k = 1}^{ \infty }$ in $A$ and a strictly increasing sequence of positive integers $\{ k(m) \}_{ m = 1}^{ \infty }$ such that
\begin{enumerate}
\item if $k(m) \leq s < k(m+1)$, then $\norm{ u_{s} a - a u_{s} } < \frac{1}{2^{ m } }$ for all $a$ in $\mc{F}_{m}$ and 
\item $[ u_{k} ]  = g_{k}$ in $U(A)/U_{0}(A)$ for all $k$.
\end{enumerate}

For each $m$, let $H_{m}$ be the finitely generated subgroup of $K_{0} ( A )$ given in Definition \ref{def:propertyC} corresponding to $\mc{F}_{m}$ and $\frac{1}{2^{m}}$.  Since $\{ G_{n} \}_{ n = 1}^{ \infty }$ is an increasing sequence of finitely generated subgroup of $K_{0} ( A )$ whose union is $K_{0} (A )$, there exists a strictly increasing sequence $\{ n(m) \}_{ m = 1}^{ \infty }$ of positive integers such that $H_{m}$ is a subgroup of $G_{ n(m) }$.  Since $\{ \check{g}_{ k } \}_{ k = 1}^{ \infty }$ converges to zero in $\invlim \frac{ K_{1} ( A ) }{ H_{ [ 1_{A} ] } ( G_{n} , K_{1} ( A ) ) }$, there exists a strictly increase sequence of positive integers $\{ k(m) \}_{ m = 1}^{ \infty }$ such that for all $k \geq k(m)$, we have that $g_{k}$ is an element of $H_{ [ 1_{A} ] } ( G_{ n(m) } , K_{1} ( A ) )$ which is a subgroup of $H_{ [ 1_{A} ] } ( H_{m} , K_{1} ( A ) )$.  By the choice $H_{m}$, for each $m$ and for $s$ in $\N$ with $k(m) \leq s < k(m+1)$, there exists a unitary $w_{s}$ in $A$ such that 
\begin{enumerate}
\item $\norm{ w_{s} a - a w_{s} } < \frac{1}{2^{m}}$ for all $a$ in $\mc{F}_{m}$ and 
\item $[ w_{s} ] = g_{s}$ in $U (A) / U_{0} ( A )$.
\end{enumerate}
For $k \geq k(1)$, set $u_{k} = w_{k}$.  For $k < k(1)$, set $u_{k}$ be any lifting of $g_{k}$.  Then, 
\begin{enumerate}
\item if $k(m) \leq s < k(m+1)$, then $\norm{ u_{s} a - a u_{s} } < \frac{1}{2^{ m }}$ for all $a$ in $\mc{F}_{m}$ and 
\item $[ u_{k} ]  = g_{k}$ in $U(A) / U_{0} (A)$ for all $k$.
\end{enumerate}
Since $\bigcup_{ m = 1}^{ \infty } \mc{F}_{m}$ is dense in $A$, it is easy to check that 
\begin{equation*}
\limit{k}{\infty}{ \norm{ u_{k} a - a u_{k} } } = 0
\end{equation*}
for all $a$ in $A$ and $[ u_{k} ] = g_{k}$ in $U ( A ) / U_{0} (A)$ for all $k$.  Therefore, by Lemma \ref{lem:stronginn} $\{ \widehat{ g_{k} } \}_{ k = 1}^{ \infty }$ converges to zero in $\frac{ \overline{ \mrm{Inn} } ( A ) }{ \overline{ \mrm{ Inn } }_{0} ( A ) }$.

The last part of the theorem is clear since $\invlim \frac{ K_{1} ( A ) }{ H_{ [ 1_{A} ] } ( G_{n} , K_{1} ( A ) ) }$ is totally disconnected.  
\end{proof}

\begin{theorem}\label{thm:inniso}
Let $A$ be a unital, separable \cstar-algebra satisfying Property (C).  Suppose the natural map from $U(A) / U_{0} (A)$ to $K_{1} (A)$ is an isomorphism.  Then for any increasing sequence $\{ G_{n} \}_{n = 1}^{ \infty }$ of finitely generated subgroups of $K_{0} ( A )$ whose union is $K_{0} ( A )$ and each $G_{n}$ contains $[ 1_{A} ]$, the topological groups $\frac{ \overline{ \mrm{Inn} } ( A ) }{ \overline{ \mrm{Inn} }_{0} ( A ) }$ and $\invlim \frac{ K_{1} ( A ) }{ H_{ [ 1_{A} ] } ( G_{n} , K_{1} ( A ) ) }$ are isomorphic. 

Moreover, if there exists an increasing sequence $\{ G_{n} \}_{ n = 1}^{ \infty }$ of finitely generated subgroup of $K_{0} ( A )$ whose union is $K_{0} ( A )$ with each $G_{n}$ containing $[ 1_{A} ]$ such that  $H_{ [ 1_{A} ] } ( G_{n} , K_{1} ( A ) ) = K_{1} ( A )$ for all $n$ sufficiently large, then $\overline{ \mrm{ Inn } } ( A ) = \overline{ \mrm{ Inn } }_{0} ( A )$.
\end{theorem}

\begin{proof}
Since the natural map from $U(A) / U_{0} (A)$ to $K_{1} (A)$ is an isomorphism, we have the following commutative diagram:
\begin{equation*}
\xymatrix{
\{1\} \ar[r] & U_{0} ( A ) \ar[r] \ar[d]_{ \mrm{Ad} } & U(A) \ar[r] \ar[d]_{ \mrm{Ad} } & K_{1} (A) \ar[r] \ar[d] & \{ 0 \} \\
\{ 1 \} \ar[r] & \overline{ \mrm{Inn} }_{0} ( A ) \ar[r] & \overline{ \mrm{Inn} } ( A ) \ar[r] & \frac{ \overline{ \mrm{Inn} } ( A ) }{ \overline{ \mrm{Inn} }_{0} ( A ) } \ar[r] & \{ 0 \}
}
\end{equation*} 
Therefore, the image of $K_{1} (A)$ is dense in the complete topological group $\frac{ \overline{ \mrm{Inn} } ( A ) }{ \overline{ \mrm{Inn} }_{0} ( A ) }$.  By Lemma \ref{lem:invlim}, $\invlim \frac{ K_{1} ( A ) }{ H_{ [ 1_{A} ] } ( G_{n} , K_{1} ( A ) ) }$ is the completion of the image of $K_{1} ( A )$.  Hence, we have that the map $\mu$ from $\frac{ \overline{ \mrm{Inn} } ( A ) }{ \overline{ \mrm{Inn} }_{0} ( A ) }$ to $\invlim \frac{ K_{1} ( A ) }{ H_{ [ 1_{A} ] } ( G_{n} , K_{1} ( A ) ) }$ defined in Theorem \ref{thm:surjinn} is surjective.  Therefore, by Theorem \ref{thm:inn}, $\mu$ is an isomorphism.  
\end{proof}

\begin{corollary}\label{cor:totallydisconnect}
Let $A$ be a separable, unital \cstar-algebra that satisfies Property (C).  Suppose the natural map from $U(A) / U_{0} ( A )$ to $K_{1} ( A )$ is an isomorphism.  Then
\begin{enumerate}

\item $\frac{ \overline{ \mrm{Inn} } ( A ) }{ \overline{ \mrm{Inn} }_{0} ( A ) }$ is compact if and only if for every finitely generated subgroup $G$ of $K_{0} ( A )$ containing $[ 1_{A} ]$, the group 
\begin{equation*}
\frac{ K_{1} ( A ) }{ H_{ [ 1_{A} ] } ( G , K_{1} ( A ) ) }
\end{equation*}
is a finite group.

\item $\frac{ \overline{ \mrm{Inn} } ( A ) }{ \overline{ \mrm{Inn} }_{0} ( A ) }$ is locally compact if and only if there exists a finitely generated subgroup $G_{0}$ of $K_{0} ( A )$ containing $[ 1_{A} ]$ such that for every finitely generated subgroup $G$ of $K_{0} ( A )$ containing $G_{0}$, the kernel of the surjective group homomorphism 
\begin{equation*}
\frac{ K_{1} ( A ) }{ H_{ [ 1_{A} ] } ( G , K_{1} ( A ) ) } \longrightarrow \frac{ K_{1} ( A ) }{ H_{ [ 1_{A} ] } ( G_{0} , K_{1} ( A ) ) }
\end{equation*}
is finite.

\item $\frac{ \overline{ \mrm{ Inn } } ( A ) }{ \overline{ \mrm{Inn} }_{0} ( A ) }$ is discrete if and only if there exists a finitely generated subgroup $G_{0}$ for $K_{0} ( A )$ containing $[ 1_{A} ]$ such that for every finitely generated subgroup $G$ of $K_{0} ( A )$ containing $G_{0}$, we have that 
\begin{equation*}
H_{ [ 1_{A} ] } ( G_{0} , K_{1} ( A ) ) = H_{ [ 1_{A} ] } ( G , K_{1} ( A ) ).
\end{equation*}
Consequently, $\frac{ \overline{ \mrm{ Inn } } ( A ) }{ \overline{ \mrm{Inn} }_{0} ( A ) }$ is isomorphic to $\frac{ K_{1} ( A ) } { H_{ [ 1_{A} ] } ( G_{0} , K_{1} ( A ) ) }$.  Also, if $\frac{ \overline{ \mrm{ Inn } } ( A ) }{ \overline{ \mrm{Inn} }_{0} ( A ) }$ is not discrete, then it has no isoloted points.
\end{enumerate}
\end{corollary}

\begin{proof}
It is easy seen that each of statement is true for the inverse limit $\invlim \frac{ K_{1} ( A ) }{ H_{ [ 1_{A} ] } ( G_{n} , K_{1} ( A ) ) }$ since the collection of subsets $\pi_{n}^{-1} ( \{ x \} )$ for $n$ in $\N$ and $x$ in $\frac{ K_{1} ( A ) } { H_{ [ 1_{A} ] } ( G_{n} , K_{1} ( A ) ) }$ is a basis of closed open sets for its topology.  Therefore, the results of the corollary follows from Theorem \ref{thm:inniso}.
\end{proof}

\section{\cstar-algebras satisfying Property (C)}\label{cstarPropC}

The propose of this section is to prove the following theorem:

\begin{theorem}\label{thm:propC}
Let $A$ be a separable, unital, amenable, simple \cstar-algebra satisfying the Universal Coefficient Theorem.  If $A$ is either a purely infinite \cstar-algebra or a tracially AI algebra, then $A$ satisfies Property (C).
\end{theorem}

The above theorem will be a consequence of the following theorem:

\begin{theorem}\label{thm:apuniteq}
Let $A$ be a unital, separable, amenable, simple \cstar-algebra which satisfies the Universal Coefficient Theorem.  Suppose $A$ is a tracially AI algebra or $A$ is a purely infinite \cstar-algebra.  Then for every $\epsilon$ in $\R_{>0}$ and finite subset $\mc{F}$ of $A$, there exists a finitely generated subgroup $G$ of $K_{0} ( A )$ containing $[ 1_{A} ]$ such that the following holds:  for every unitary $u$ in $A$ with $[ u ]$ an element of $H_{ [ 1_{A} ] } ( G , K_{1} ( A ) )$, there exist a unital, contractive, completely positive, linear map $\psi$ from $A \otimes C(S^{1})$ to $A$ and a unitary $w$ in $A$ such that $\Phi$ is $\iota ( \mc{F} ) \cup \{ 1_{A} \otimes z , 1_{A} \otimes z^{*} \}$-$\epsilon$-multiplicative, $K_{1}( \Phi ) ( [ 1_{A} \otimes z ] ) = [ u ] $ in $K_{1} ( A )$, and
\begin{equation*}
\norm{ w ( \Phi \circ \iota ) ( a ) w^{*} - a  } < \epsilon
\end{equation*}
for all $a$ in $\mc{F}$.
\end{theorem}

Let us use Theorem \ref{thm:apuniteq} to prove Theorem \ref{thm:propC}.  We will first need the following well-known lemma (see Lemma 4.1.1 of \cite{booklin}).

\begin{lemma}\label{lem:almostunit}
Let $\epsilon$ be an element of $\R_{ > 0}$.  Then, there exists $\delta$ in $\R_{>0}$ satisfying the following:  For any unital \cstar-algebra $A$ and $x$ in $A$, if $\norm{ x^{*} x - 1_{A} } < \delta$ and $\norm{ x x^{*} - 1_{A} } < \delta$, then there exists a unitary $u$ in $A$ such that 
\begin{equation*}
\norm{ u - x } < \epsilon.
\end{equation*}
\end{lemma}

Note that if $\epsilon$ is strictly less than $\frac{1}{2}$ and $w$ is another unitary in $A$ such that 
\begin{equation*}
\norm{ w - x } < \epsilon,
\end{equation*}
then $ [ w ] = [ u ]$ in $K_{1} (A)$.  Hence, in Theorem \ref{thm:apuniteq} with $\epsilon$ stricitly less than $\frac{1}{2}$, we have that $K_{1} ( [ 1_{A} \otimes z ] )$ is a well-define element of $K_{1} ( A )$.

\bigskip

\noindent\emph{Proof of Theorem \ref{thm:propC}.}\ Let $\epsilon$ be in $\R_{ > 0}$ and $\mc{F}$ be a finite subset of $A$.  Note that we may assume that $\epsilon$ is strictly less than $\frac{1}{2}$ and $\mc{F}$ is a finite subset of the unit ball of $A$.  Let $\delta$ in $\R_{ > 0}$ be the positive number provided by Lemma \ref{lem:almostunit} corresponding to $\frac{ \epsilon }{ 6 }$.  Note that we may assume that $\delta$ is strictly less than $\frac{ \epsilon }{ 6 }$.  

Let $G$ be the finitely generated subgroup of $K_{0} ( A )$ containing $[ 1_{A} ] $ given by Theorem \ref{thm:apuniteq} which corresponds to $A$, the finite set $\mc{F}$, and the positive number $\delta$.  Hence, if $u$ is a unitary in $A$ with $[ u ]$ an element of $H_{ [ 1_{A} ] } ( G , K_{1} ( A ) )$, then by Theorem \ref{thm:apuniteq} there exist a unital, completely positive, linear map $\Phi$ from $A \otimes C(S^{1})$ to $A$, and a unitary $w$ in $A$ such that
\begin{enumerate}
\item $\Phi$ is $\iota ( \mc{F} ) \cup \{ 1_{A} \otimes z , 1_{A} \otimes z^{*} \}$-$\delta$-multiplicative;
\item $K_{1} ( \Phi ) ( [ 1_{A} \otimes z ] ) = [ u ]$ in $K_{1} ( A )$;
\item $\norm{ w ( \Phi \circ \iota )( a ) w^{*} - a } < \delta$ for all $a$ in $\mc{F}$.
\end{enumerate}

By the choice of $\delta$ and by Lemma \ref{lem:almostunit}, there exists a unitary $v_{0}$ in $A$ such that 
\begin{equation*}
\norm{ \Phi ( 1_{A} \otimes z ) - v_{0} } < \frac{ \epsilon }{ 6 }
\end{equation*}
Set $v = w v_{0} w^{*}$.  Then $[ v ] = [ v_{0} ] = K_{1}( \Phi ) ( [ 1_{A} \otimes z ] ) = [ u ] $ in $K_{1} ( A )$ and hence, by Proposition \ref{prop:k1injTAI} if $A$ is a tracially AI algebra and by Theorem 1.9 of \cite{kthypureinf} if $A$ is a purely infinite simple \cstar-algebra, $ [ v ] = [ u ]$ in $U(A) / U_{0} ( A )$.  A computation shows that
\begin{equation*}
\norm{ v a - a v } < \epsilon
\end{equation*} 
for all $a$ in $\mc{F}$.  

\bigskip

In order to prove Theorem \ref{thm:apuniteq}, we first develop some notation that will be used through the rest of this section.  Let $A$ be a unital \cstar-algebra and let $C_{n}$ be the mapping cone of the degree $n$ map on $C_{0}(\R)$.  The set of all projections and unitaries in ${ \displaystyle \bigcup_{m,n=1}^{\infty} \msf{M}_{m}( \unitize{ A \otimes C_{n}} ) }$ will be denoted by $\proj(A)$.  Denote the total $K$-theory of $B$ defined in \cite{approxhom} by $\totalk(B)$.  For any finite subset $\mc{P}$ of $\proj(A)$, there exist a finite subset $\mc{G}(\mc{P})$ of $A$ and $\delta(\mc{P})$ is $\R_{ > 0 }$ such that if $B$ is any \cstar-algebra and $L$ is a contractive, completely positive, linear map from $A$ to $B$ which is $\mc{G}(\mc{P})$-$\delta(\mc{P})$-multiplicative, then $L$ defines a map $\totalk( L ) \vert_{ \mc{P} }$ from $\overline{ \mc{P} }$ to $\underline{K}(B)$, where $\overline{\mc{P}}$ is the image of $\mc{P}$ in $\totalk(A)$.  By enlarging $\mc{G}(\mc{P})$ and chosing a smaller $\delta ( \mc{P} )$, if necessary, $\totalk( L ) \vert_{ \mc{P} }$ is defined on the subgroup generated by $\overline{ \mc{P} }$.  

Suppose $A$ satisfies the Universal Coefficient Theorem.  Then, by Theorem 1.4 of \cite{multcoeff}, $\kl(A,B)$ is naturally isomorphic to $\Hom_{\Lambda}(\totalk(A), \totalk(B))$ and under this isomorphism $\kl ( \varphi )$ is precisely $\totalk( \varphi )$ for any $*$-homomorphism $\varphi$ from $A$ and $B$.  In the sequel, we will make this identification without further mention.   

We will also need the following lemma.

\begin{lemma}\label{lem:exthom}
Let $A$ be a unital, separable, amenable \cstar-algebra satisfying the Universal Coefficient Theorem and let $B$ be a unital \cstar-algebra.  Suppose $\varphi$ is an element of $\Hom_{ \Lambda } ( \totalk(A) , \totalk(B) )$ such that $\varphi( [1_{A}] ) = [ 1_{B} ]$ and $\varphi \vert_{ K_{0} ( A ) }$ is a positive homomorphism and suppose $\zeta$ is an element of $H_{ [ \unit{A} ] }(K_{0}(A), K_{1}(B))$.  Then, there exists $\alpha$ in $\Hom_{ \Lambda } ( \totalk( A \otimes C(S^1) ) , \totalk( B ) )$ such that
\begin{enumerate}
\item  $\alpha \vert_{ K_{0} ( A \otimes C(S^{1} ) ) }$ is positive;
\item  $\alpha \circ \totalk( \iota ) = \totalk ( \varphi )$; and
\item  $\alpha ( [ \unit{A} \otimes z ] ) = \zeta$.
\end{enumerate}
\end{lemma}

\begin{proof}
Throughout the proof $a \times b$ will denote the Kasparov product of $a$ and $b$.  Note that the sequence
\begin{equation*}
\xymatrix{ 0 \ar[r] & C_{0}((0,1) , A) \ar[r]  & A \otimes C( S^{1} ) \ar[r]^(.65){\pi}  & A \ar[r] & 0}
\end{equation*}  
is a split exact sequence, where the splitting map is given by $\iota ( a ) =  a \otimes 1_{ C( S^{1} ) }$.  Hence, $K_{0} ( A \otimes C(S^{1}) ) \cong K_{0} (A) \oplus K_{1} (A)$ and $K_{1} ( A \otimes C(S^{1}) ) \cong K_{1} (A) \oplus K_{0} (A)$.  Moreover, under these isomorphisms, $K_{0}(\iota)$ from $K_{0} (A)$ to $K_{0} (A) \oplus K_{1} (A)$ may be written as $K_{0}(\iota) (x) = ( x , 0)$ and $K_{1}(\iota)$ from $K_{1}(A)$ to $K_{1}(A) \oplus K_{0}(A)$ may be written as $K_{1}(\iota) (x) = (x ,0)$.  

Define $\gamma_{0}$ from $K_{0} (A \otimes C(S^{1}))$ to $K_{0}(B)$ by $\gamma_{0} = \varphi \circ K_{0}(\pi)$.  Since $\zeta$ is an element of $H_{ [ \unit{A} ] } (K_{0} (A), K_{1} (B))$, there exists a homomorphism $\alpha$ from $K_{0}(A)$ to $K_{1}(B)$ such that $\alpha( [ 1_{A} ]) = \zeta$.  Define $\gamma_{1}$ from $K_{1} (A \otimes C(S^{1}) )$ to $K_{1}(B)$ by $\gamma_{1} ((x,y)) = \varphi (x) + \alpha(y)$.

Since $A$ satisfies the Universal Coefficient Theorem, $A \otimes C(S^{1})$ satisfies the Universal Coefficient Theorem.  Hence, there exists $\xi$ in $\kl(A \otimes C(S^{1}), B)$ such that $\Gamma(\xi)_{i} = \gamma_{i}$.  Therefore, $\varphi - \kl( \iota ) \times \xi$ is an element of $\ext(K_{*}(A), K_{*+1}(B))$.  Suppose $\varphi - \kl ( \iota ) \times \xi$ is represented by the following exact sequences:
\begin{align*}
\stex{K_{1}(B)}{}{K_{0}(E)}{}{K_{0}(A)}\\
\stex{K_{0}(B)}{}{K_{1}(E)}{}{K_{1}(A)}.
\end{align*}
Let $G_{0} = K_{0}(E) \oplus K_{1} (A)$ and $G_{1} = K_{1}(E) \oplus K_{0}(A)$.  Then the sequence
\begin{equation*}
\stex{K_{i + 1}(B) }{}{G_{i}}{}{K_{i}(A) \oplus K_{i+1} (A)} 
\end{equation*}
is exact for $i = 0 ,1$, which represents an element $\delta$ of $\ext( K_{*} (A \otimes C(S^{1})) , K_{* + 1} (B) )$ such that $\kl ( \iota ) \times \delta = \varphi - \kl( \iota ) \times \xi$.  Let $\alpha = \xi + \delta$.  It is easy to check that $ \varphi = \kl (\iota ) \times \alpha $ and $\Gamma(\alpha)_{1}( [ 1_{A} \otimes z ] ) = \Gamma(\xi)_{1}([ 1_{A} \otimes z])  = \gamma_{1}(0,1)  = \zeta$.  So $\phi = \alpha \circ \totalk( \iota )$ and $\phi (  [ 1_{A} \otimes z ] ) = \zeta$.  Note that $\Gamma( \alpha )_{0} = \Gamma ( \xi )_{0} + \Gamma ( \delta )_{0} = \gamma_{0} =  \varphi \circ K_{0} ( \pi )$.  Thus, $\alpha \vert_{ K_{0} ( A \otimes C(S^{1}) ) }$ is a positive homomorphism.
\end{proof}

\subsection{Proof of Theorem \ref{thm:apuniteq}: The Purely Infinite Case} Let us recall the following result which is an easy consequence of the results of Kirchberg \cite{kirchpureinf}, Phillips \cite{phillipspureinf}, and Elliott and R{\o}rdam \cite{elliottrordam} and \cite{defKL}.  

\begin{theorem}\label{thm:pureinfdecomp}
If $A$ is a simple, amenable, separable, unital purely infinite \cstar-algebra satisfying the Universal Coefficient Theorem, then $A$ is isomorphic to $\dirlim ( A_{n} , \varphi_{n,n+1} )$, where each $A_{n}$ is a unital, separable, amenable, purely infinite simple \cstar-algebra with finitely generated $K$-theory and satisfies the Universal Coefficient Theorem and each $\varphi_{n, n+1}$ is unital.  
\end{theorem}
  
\noindent\emph{Proof of Theorem \ref{thm:apuniteq}: The Purely Infinite Case.}  Let $\epsilon$ an element of $\R_{ > 0}$ and $\mc{F}$ a finite subset of $A$.  By Theorem \ref{thm:pureinfdecomp}, $A$ is isomorphic to $\dirlim ( A_{n} , \phi_{ n , n+1} )$ where each $A_{n}$ is a unital, separable, amenable, purely infinite simple \cstar-algebra with finitely generated $K$-theory and satisfies the Universal Coefficient Theorem and each $\varphi_{n, n+1}$ is unital.  Hence, there exists a sub-\cstar-algebra $B$ of $A$ and a finite subset $\mc{F}_{B}$ of $B$ such that 
\begin{enumerate}
\item[(i)]  $B$ is a separable, unital, amenable, purely infinite simple \cstar-algebra satisfying the Universal Coefficient Theorem;
\item[(ii)]  $1_{B} = 1_{A}$;
\item[(iii)]  $K_{*} ( B )$ is finitely generated; and
\item[(iv)]  every element of $\mc{F}$ is within $\frac{ \epsilon }{ 7 }$ of an element of $\mc{F}_{B}$.
\end{enumerate} 

Denote the inclusion of $B$ into $A$ by $\psi$.  Set $G_{0} =  K_{0} ( \psi ) ( K_{0} ( B ) )$.  Suppose $u$ is a unitary in $A$ such that $[ u ]$ is an element of $H_{ [ 1_{A} ] } ( G_{0}, K_{1} ( A ) )$.  Since $\psi( 1_{B} ) = 1_{B} = 1_{A}$, we have that $[ u ]$ is an element of $H_{ [ 1_{B} ] } ( K_{0} (B) , K_{1} ( A ) )$.   By Lemma \ref{lem:exthom}, there exists $\alpha$ in $\Hom_{ \Lambda } ( \totalk( B \otimes C(S^{1}) ) , \totalk( A ) )$ such that 
\begin{enumerate}
\item $\alpha \circ \totalk ( \iota ) = \totalk ( \psi )$ and 
\item $\alpha ( [ 1_{ B } \otimes z ] ) = [ u ]$.
\end{enumerate}     
By Theorem 6.7 of \cite{sepBDF}, there exists a unital injective $*$-homomorphism $\beta$ from $B \otimes C(S^{1})$ to $A$ such that $\totalk ( \beta ) = \alpha$.  By Theorem 6.7 of \cite{sepBDF}, there exists a unitary $w$ in $A$ such that 
\begin{equation*}
\norm{ w ( \beta \circ \iota ) ( b ) w^{*} - \psi ( b ) } < \frac{ 4 \epsilon }{ 7 }
\end{equation*}
for all $b$ in $\mc{F}_{B}$.  Since $\psi ( b ) = b$ for all $b$ in $B$, we have that 
\begin{equation*}
\norm{ w ( \beta \circ \iota ) ( b ) w^{*} - b } < \frac{ 4 \epsilon }{ 7 }
\end{equation*} 
for all $b$ in $\mc{F}_{B}$.

Since $B \otimes C(S^{1})$ is an amenable sub-\cstar-algebra of $A \otimes C(S^{1})$, there exists a unital, contractive, completely positive, linear map $\Phi$ from $A \otimes C(S^{1})$ to $A$ such that 
\begin{equation*}
\norm{ \Phi ( x ) - \beta ( x ) } < \frac{ \epsilon }{ 7 }
\end{equation*}
for all $x$ in the set
\begin{equation*}
\iota ( \mc{F}_{B} ) \cup \{ 1_{A} \otimes z , 1_{A} \otimes z^{*} \} \cup \setof{ x y }{ x , y \in \iota ( \mc{F}_{B} ) \cup \{ 1_{A} \otimes z , 1_{A} \otimes z^{*} \} }.
\end{equation*}
A computation shows that $\Phi$ is $\mc{F} \cup \{ 1_{A} \otimes z , 1_{A} \otimes z^{*} \}$-$\epsilon$-multiplicative, $K_{1} ( \Phi )( [ 1_{A} \otimes z ] ) = [ u ]$ in $K_{1} ( A )$, and 
\begin{equation*}
\norm{ w ( \Phi \circ \iota )( a ) w^{*} - a } < \epsilon 
\end{equation*}
for all $a$ in $\mc{F}$.

\subsection{Proof of Theorem \ref{thm:apuniteq}: The tracially AI case}

Let us recall the definition of a tracially AI algebra due to Lin \cite{linTR1}.

\begin{definition}\label{def:TAI}
Let $\mc{I}$ be the class of all unital \cstar-algebras with the form $\oplus_{i=1}^{n} B_{i}$, where each $B_{i}$ is isomorphic to $\msf{M}_{k}$ or $\msf{M}_{k} ( C( [0,1] ) )$ for some integer $k$.  

A unital \cstar-algebra $A$ is said to be \emph{tracially AI} if for any finite subset $\mc{F}$ of $A$
containing a nonzero element $b$, $\epsilon$ in $\R_{ > 0 }$, a positive integer $n$, and any positive element $a$ of $A$ which is not contained in any proper ideal of $A$, there exists a nonzero projection $p$ in $A$ and a sub-\cstar-algebra $I$ of $A$ with $I$ in $\mc{I}$ and $1_{I} = p$ satisfying the following properties: 
\begin{enumerate}
\item $\norm{ p x - x p } < \epsilon$ for all $x$ in $\mc{F}$;

\item for every $x$ in $\mc{F}$, we have that $p x p$ is to within $\epsilon$ of an element of $I$ and $\norm{ pbp } \geq \norm{ b } - \epsilon$; and 

\item $n(1_{A} - p)$ is Murray-von Neumann equivalent to a sub-projection of $p$ and $1_{A} - p$ is Murray-von Neumann equivalent to a projection in $\overline{ a A a }$.
\end{enumerate}

A non-unital \cstar-algebra $A$ is said to be tracially AI if $\widetilde{ A}$ is a tracially AI algebra.

By Theorem 4.10 of \cite{linTR1}, if $A$ is simple, condition (3) can replaced by

(3a) $1_{A}-p$ is unitarily equivalent to a projection in $eAe$ for any previously given nonzero projection $e$ in $A$

If $A$ has the Fundamental Comparability Property (see \cite{BlackComp}),
condition (3) can be replaced by

(3b) $\tau(1_{A} - p) < \sigma$ for any prescribed $\sigma$ in $\R_{ > 0}$ and for
all normalized quasi-traces of $A.$

A \cstar-algebra $A$ is said to be a \emph{tracially AF algebra} if the class of sub-\cstar-algebras $\mc{I}$ in the definition of a tracially AI algebra is replaced by the class of all finite dimensional \cstar-algebras.
\end{definition}

We will first show that the canonical map from $U(A) / U_{0} (A)$ to $K_{1} (A)$ is an isomorphism from simple, unital, tracially AI algebras.  We thank Zhuang Niu for this argument.

\begin{lemma}\label{lem:contunitary}
Let $u$ and $v$ be unitaries in a unital \cstar-algebra $A$.  If $\norm{ u - v } < 2$ and $p$ is a projection in $A$ such that $p$ commutes with $u$ and $v$, then there exists a continuous path of unitaries $w(t)$ in $A$ such that $w(0) = v$, $w(1) = u$, and $p$ commutes with $w(t)$ for all $t$ in $[ 0 , 1 ]$.  Moreover, if $p u p$ and $p v p$ are in a unital sub-\cstar-algebra $B$ of $p A p$, then $w(t)$ can be chosen such that $p w(t) p$ is a continuous path of unitaries in $B$.
\end{lemma}

\begin{proof}
Since $\norm{ u - v } < 2$, we have that $-1$ is not in the spectrum of $v^{*} u$.  Therefore, there exists a real valued continuous function $\phi$ on the spectrum of $v^{*} u$ such that $z = \exp ( i \phi (z) )$ for all $z$ in the spectrum of $v^{*}u$.  Set $h = \phi ( v^{*} u )$ which is a self-adjoint element of $A$.  Note that $v^{*} u = \exp ( i h )$.  

Now set $w_{0}( t ) = \exp( i t h )$.  Then $w_{0} ( t )$ is a continuous path of unitaries in $A$ such that $w_{0} ( 0 ) = 1_{A}$ and $w_{0} ( 1 ) = v^{*} u$.  Since $p$ commutes with $u$ and $v$, we have that $p$ commutes with $h$.  Hence, $p$ commutes with $w_{0} ( t )$ for all $t$.  Set $w(t) = v w_{0} (t)$.  Then $w(t)$ is a continuous path of unitaries in $A$ such that $w(0) = v$, $w(1) = u$, and $p$ commutes with $w(t)$ for all $t$ in $[ 0 , 1 ]$.  

Suppose $ p u p$ and $p v p $ are in a unital sub-\cstar-algebra $B$ of $p A p$.  Then one has that 
\begin{equation*}
p  h p = \phi ( (p v^{*} p) ( p u p ) )
\end{equation*} 
which is an element of $B$.  Hence, 
\begin{equation*}
p w(t) p = p v w_{0} (t) p = p v p \exp ( i t p h p )
\end{equation*}
which is an element of $B$.
\end{proof}

\begin{proposition}\label{prop:k1injTAI}
Let $A$ be a simple tracially AI algebra.  Then the canonical map from $U(A) / U_{0} (A)$ to $K_{1} (A)$ is an isomorphism. 
\end{proposition}

\begin{proof}
Since $A$ is a simple tracially AI algebra, by Theorem 4.5 of \cite{linTR1} $A$ has stable rank one.  Hence, the canonical map from $U(A) / U_{0} (A)$ to $K_{1} (A)$ is surjective (see Lemma 3.1.10 of \cite{booklin}).  To show this map is injective, it is enough to show that for any unitary $u$ in $A$ if 
\begin{equation*}
\mrm{ diag } ( u , 1_{A} , \dots , 1_{A} )
\end{equation*}
is path connected to 
\begin{equation*}
\mrm{ diag } ( 1_{A} , \dots, 1_{A} )
\end{equation*}
by a continuous path of unitaries in $\msf{M}_{ n + 1} (A)$, then $u$ is connected to $1_{A}$ by a continuous path of unitaries in $A$.  

Suppose we have a continuous path of unitaries $W(t)$ in $\msf{M}_{ n + 1} (A)$ connecting $\mrm{ diag } ( u , 1_{A} , \dots, 1_{A} )$ to $\diag ( 1_{A} , \dots, 1_{A} )$.  Then there exists a partition $0 = t_{1} < t_{2} < \cdots < t_{s} = 1$ of $[ 0 , 1 ]$ such that 
\begin{equation*}
\norm{  W ( t_{k} ) - W ( t_{ k+1} ) } < 2
\end{equation*}
for all $0 \leq k \leq s - 1$.  

Note that by Theorem 3.2 of \cite{linTR1}, every nonzero hereditary sub-\cstar-algebra of $A$ contains a nonzero projection.  Therefore, there are $n + 1$ mutually orthogonal projections $q_{1}, \dots, q_{n+1}$ in $A$ which is Murray-von Neumann equivalent to each other (see Lemma 3.5.7 of \cite{booklin}).  Now since $A$ is a simple, tracially AI algebra we may assume that there exists a sub-\cstar-algebra $C$ in $\mc{I}$ with $ 0 \neq p = 1_{ C } $ such that 
\begin{equation*}
W( t_{k} )_{ i , j } = W_{1} ( t_{k} )_{ i , j } + W_{2} ( t_{k} )_{ i , j }
\end{equation*} 
with $W_{1} ( t_{k} )_{ i , j}$ an element of $( 1_{A} - p ) A ( 1_{A} - p )$ and $W_{2}( t_{k} )_{ i , j }$ an element of $C$ for each $0 \leq k \leq s$ and $1 \leq i, j \leq n+1 $, and $1_{A} - p$ is Murray-von Neumann equivalent to a subprojection of $q_{1}$.  Moreover, one can assume that $( W_{2} ( t_{k} )_{ i , j } )$ is a unitary in $\msf{M}_{ n +1 } ( ( 1_{A} - p ) A ( 1_{A} - p ) )$ and $( W_{2} ( t_{k} )_{ i, j } )$ is a unitary in $\msf{M}_{n + 1} ( C )$.  By Lemma \ref{lem:contunitary}, there exists a continuous path of unitaries $w_{k} ( t )$ between $W( t_{k} )$ and $W( t_{k+1} )$ for each $0 \leq k \leq s - 1$ and $w_{k}$ commutes with $Q = \diag ( 1_{A} - p , 1_{A} - p , \dots , 1_{A} - p )$.  Therefore, $Q w_{k} ( t ) Q$ is a continuous path of unitaries in $\msf{M}_{ n + 1 } ( ( 1_{A} - p ) A ( 1_{A} - p ) )$ which connects $( W_{1} ( t_{k} )_{i , j } )$ and $( W_{1} ( t_{k+1} )_{ i , j } )$.  Moreover, $ ( 1_{ \msf{M}_{n+1} ( A ) } - Q ) w_{k} ( t ) ( 1_{ \msf{M}_{n+1} ( A ) } - Q )$ is a continuous path of unitaries in $\msf{M}_{ n + 1} ( C )$ which connects $( W_{2} ( t_{k} )_{ i , j } )$ and $( W_{2} ( t_{k+1} )_{ i , j } )$.  By patching these continuous paths of unitaries together, we get a continuous path of unitaries $w(t)$ such that $Q w(t) Q$ is a continuous path of unitaries in $\msf{M}_{ n + 1} ( ( 1_{A} - p ) A ( 1_{A} - p ) )$ which connects $( W_{1} ( t_{0} )_{ i , j } )$ and $( W_{1} ( t_{s} )_{ i , j } ) = Q$ and $( 1_{ \msf{M}_{n+1} ( A ) } - Q ) w(t) ( 1_{ \msf{M}_{n+1} ( A ) } - Q )$ is a continuous path of unitaries in $\msf{M}_{ n + 1} ( C )$ which connects the unitary
\begin{equation*}
( W_{2} ( t_{0} )_{ i , j } ) = \diag ( W_{2} ( t_{0} )_{1,1} , p , \dots , p )
\end{equation*}
to the unitary
\begin{equation*}
( W_{2} ( t_{s} )_{ i , j } ) = ( 1_{ \msf{M}_{n+1} ( A ) } - Q ) = \diag (  p , \dots , p ).
\end{equation*}
Therefore, the unitary $W_{2} ( t_{0} )_{1,1}$ is zero in $K_{1} (C)$.  Since the canonical map from $U(C) / U_{0} ( C )$ to $K_{1} ( C )$ is injective, we have that $W_{2} ( t_{0} )_{1,1}$ is in $U_{0} ( C )$.

Recall that $1_{A} - p$ is Murray-von Neumann equivalent to a subprojection of $q_{1}$ and the mutually orthogonal projections $q_{1}, \dots, q_{ n + 1}$ are Murray-von Neumann equivalent to each other.  Since $A$ has stable rank one, there are partial isometries $v_{1} = 1_{A} - p$, $v_{2}, \dots , v_{n+1}$ such that the source projections are equal to $1_{A} - p$ and the range projections are mutually orthogonal.  Set
\begin{equation*}
V = ( v_{1} , v_{2} , \dots v_{ n +1 } )
\end{equation*}
as an element of $\msf{M}_{ 1, n+1} (A )$.  It is easy to check that $V^{*} V = Q$.  Then 
\begin{equation*}
c(t) = V Q w(t) Q V^{*} + ( 1_{A} - V V^{*} )
\end{equation*}
is a continuous path of unitaries in $A$.  Note that 
\begin{align*}
c(0) &= V \diag ( ( 1_{A} - p ) u ( 1_{A} - p ) , 1_{A} - p , \dots , 1_{A} - p ) V^{*} + 1_{A} - V V^{*} \\
	&= ( 1_{A} - p ) u ( 1_{A} - p ) + p \\
	&= W_{1} ( t_{0} )_{1,1} + p
\end{align*}
and 
\begin{align*}
c(1) &= V \diag ( 1_{A} - p , 1_{A} - p , \dots , 1_{A} - p ) V^{*} + 1_{A} - V V^{*} \\
	&= 1_{A} - p + p \\
	&= 1_{A}.
\end{align*}
Therefore, $W_{1} ( t_{0} )_{1,1} + p$ is connected to $1_{A}$ by a continuous path of unitaries in $A$.  

Note that 
\begin{equation*}
u = W( t_{0} )_{1,1} = W_{1} ( t_{0} )_{1,1} + W_{2} ( t_{0} )_{1,1}.
\end{equation*}
Since $W_{2} ( t_{0} )_{1,1}$ is connected to $p = 1_{C}$ by a continuous path of unitaries in $C$, one has that $u$ is connected to $c(0) = W_{1} ( t_{0} )_{1,1} + p$ by a continuous path of unitaries in $A$, and hence $u$ is connected to $1_{A}$ by a continuous path of unitaries in $A$.
\end{proof}

\begin{lemma}\label{lem:pertunital}
Let $A$ be a unital \cstar-algebra and let $\mc{F}$ be a finite subset of $A$ containing the identity of $A$ and $\mc{P}$ a finite subset of $\proj ( A )$.  For every $\epsilon$ in $\R_{ > 0}$, there exists $\delta$ in $\R_{>0}$ such that the following holds:  If $B$ is a \cstar-algebra and $\phi$ from $A$ to $B$ is $\mc{F}$-$\delta$-multiplicative, contractive, completely positive, linear map with $\totalk( \phi ) \vert_{ \mc{P} }$ well-defined, then there exists a $\mc{F}$-$\epsilon$-multiplicative, contractive, completely positive, linear map $\psi$ from $A$ to $B$ such that $\psi ( 1_{A} )$ is a projection in $\mc{B}$,	
\begin{equation*}
\norm{ \phi ( a ) - \psi ( a ) } < \epsilon
\end{equation*}
for all $a$ in $\mc{F}$ and $\totalk ( \phi ) \vert_{ \mc{P} } = \totalk ( \psi ) \vert_{ \mc{P} }$.
\end{lemma}

\begin{lemma}\label{lem:liftinghom}
Let $A$ be a separable, unital, amenable, simple \cstar-algebra satisfying the Universal Coefficient Theorem.  Suppose $A$ is a tracially AI algebra.  Then for every $\epsilon$ in $\R_{ > 0}$, finite subset $\mc{P}$ of $\proj ( A )$, and finite subset $\mc{F}$ of $A \otimes C( S^{1} )$, there exists a finitely generated subgroup $G$ of $K_{0} ( A )$ containing the class of the identity of $A$ such that the follow holds: for every unitary $u$ in $A$ with $[ u ] $ an element of $ H_{ [ 1_{A}] } ( G , K_{1} ( A ) )$, there exists a contractive, completely positive, linear map $\psi$ from $A \otimes C( S^{1} )$ to $A$ such that 
\begin{enumerate}
\item $\psi$ is $\mc{F}$-$\epsilon$-multiplicative;
\item $\totalk ( \psi \circ \iota ) \vert_{ \mc{P} } = \totalk ( \id_{A} ) \vert_{ \mc{P} } $; and
\item $\totalk ( \psi ) ( [ 1_{A} \otimes z ] ) = [ u ]$.
\end{enumerate}  
\end{lemma}

\begin{proof}
Note that without loss of generality, we may assume that $\mc{F}$ is contained in the unit ball of $A \otimes C( S^{1} )$.  We also may assume that $1_{A} \otimes z$ and $1_{A} \otimes z^{*}$ are in $\mc{F}$ and $\epsilon$ is strictly less than one.  Set $\mc{Q} = \iota ( \mc{P} ) \cup \{ 1_{A} \otimes z, 1_{A} \otimes 1_{ C( S^{1} ) } \}$.   

Suppose $A$ is finite dimensional.  Then $A = \msf{M}_{ n } ( \C )$.  Take $G = K_{0} ( A )$.  Note that $K_{1} ( A ) = 0$.  Thus, if $ [ u ]$ is an element of $H_{ [ 1_{A} ] } ( K_{0} ( A ) , K_{1} ( A ) ) = 0$, then $ [ u ] = 0$.  Therefore, we can take $\psi$ to be the natural projection from $A \otimes C(S^{1})$ to $A$.

Suppose $A$ is infinite dimensional.  By Theorem 10.9 and Theorem 10.1 of \cite{linTR1}, $A$ is isomorphic to $\dirlim ( A_{n} , \phi_{ n, n+1} )$, where $A_{n} = \bigoplus_{ i = 1 }^{ k(n) } P_{ [ n, i ] } \msf{M}_{ [ n , i ] } ( C( X_{ [n,i] }  ) ) P_{ [n,i] }$ with each $X_{ [n,i] }$ a connected finite CW-complex and $\phi_{n,n+1}$ is unital and injective.  So we may assume that $A$ is this direct limit decomposition.  By \cite{TR0} and \cite{linTR1}, there exists a unital, separable, simple AH algebra $B$ such that $B$ is a tracially AF algebra and
\begin{equation*}
( K_{0} ( A ) , K_{0} ( A )_{ + } , [ 1_{A} ] , K_{1} ( A ) ) \cong ( K_{0} ( B ) , K_{0} ( B )_{ + } , [ 1_{B} ] , K_{1} ( B ) ) .
\end{equation*}
Denote the above isomorphism by $\beta$.  Since $A$ satisfies the Universal Coefficient Theorem, $\beta$ lifts to an isomorphism from $\totalk ( A )$ to $\totalk ( B )$.  By an abuse of notation, we denote this lifting by $\beta$.  

Denote the natural embedding of $A_{n}$ into $A_{n} \otimes C(S^{1})$ by $\iota_{ A_{n} }$.  Since $ A \otimes C(S^{1}) = \dirlim ( A_{n} \otimes C(S^{1}) , \phi_{ n , n+1} \otimes \id_{ C(S^{1} ) } )$ and since $\totalk ( \cdot )$ is a continuous functor with respect to direct limits, there exist $n_{0}$ in $\N$, finite subset $\mc{P}_{ n_{0} } $ of $\proj ( A_{ n_{0} } )$, finite subset $\mc{Q}_{n_{0}}$ of $\proj ( A_{n_{0}} \otimes C(S^{1}) )$ containing $\iota_{ A_{ n_{0} } } ( \mc{P}_{ n_{0} } )$, and a finite subset $F_{ n_{0} }$ of $A_{n_{0}} \otimes C(S^{1})$ such that 
\begin{enumerate}
\item for each $p$ in $\mc{P}$, there exists $e_{p}$ in $\mc{P}_{ n_{0} }$ such that $[ p ] = \totalk ( \phi_{ n_{0} , \infty } ) ( [ e_{p} ] )$;
\item for each $p$ in $\mc{Q}$, there exists $q_{p}$ in $\mc{Q}_{ n_{0} }$ such that $[ p ] = \totalk ( \phi_{ n_{0} , \infty } \otimes \id_{ C(S^{1} ) } )( [ q_{p} ] )$;
\item every element of $\mc{F}$ is within $\frac{ \epsilon }{ 20 }$ to an element of $( \phi_{ n_{0} , \infty } \otimes \id_{ C(S^{1} ) } )( \mc{F}_{ n_{0} } )$.  
\end{enumerate}

Set $G = K_{0} ( \phi_{ n_{0} , \infty } ) ( K_{0} ( A_{ n_{0} } ) )$.  Suppose $[ u ]$ is an element of $H_{ [ 1_{A} ] } ( G , K_{1} ( A ) )$.  Then $\beta ( [ u ] )$ is an element of $H_{ [ 1_{ A_{ n_{0} } } ] } ( K_{0} ( A_{ n_{0} } ) , K_{1} ( B ) )$.  Therefore, by Lemma \ref{lem:exthom} there exists $\alpha$ in $\Hom_{ \Lambda } ( \totalk ( A_{ n_{0} } \otimes C( S^{1} ) ) , \totalk ( B ) )$ such that
\begin{enumerate}
\item $\alpha \vert_{ K_{0} ( A_{ n_{0} } \otimes C(S^{1}) ) }$ is positive;
\item $\alpha \circ \totalk( \iota_{ A_{ n_{0} } } ) = \beta \circ \totalk ( \phi_{ n_{0} , \infty } )$; and
\item $\alpha ( [ 1_{ A_{n_{0}} }  \otimes z ] ) = \beta ( [  u ] ) $ in $K_{1} ( B )$.
\end{enumerate}
By Theorem 6.2.9 of \cite{booklin}, there exists a sequence of contractive, completely positive, linear maps $\{ L_{ n_{0} , k} \}_{ k = 1 }^{\infty}$ from $A_{ n_{0} } \otimes C( S^{1} )$ to $B$ such that 
\begin{enumerate}
\item $\limit{ k }{ \infty } { \norm{ L_{ n_{0} , k } ( x y ) - L_{ n_{0} , k } ( x ) L_{ n_{0} , k } ( y ) } } = 0$ for all $x$ and $y$ in $A_{ n_{0} } \otimes C( S^{1} )$ and
\item $\totalk ( L_{ n_{0} , k } ) \vert_{ \mc{Q}_{ n_{0} } } = \alpha \vert_{ \mc{Q}_{ n_{0} } }$ for all $k$.
\end{enumerate}
Since $\totalk ( L_{ n_{0} , k } ) ( [ 1_{ A_{n_{0}} } \otimes 1_{ C(S^{1}) } ] ) = \alpha ( [ 1_{ A_{n_{0}} } \otimes 1_{ C(S^{1}) } ] ) = [ 1_{B} ]$, by a small perturbation we may assume that $L_{ n_{0} , k }$ is unital for all $k$.  Since $A_{ n_{0} } \otimes C( S^{1} )$ is amenable, there exists a sequence of contractive, completely positive, linear maps $\{ \psi_{ n_{0} , k } \}_{ k = 1 }^{ \infty }$ from $A \otimes C( S^{1} )$ to $ A_{ n_{0} } \otimes C( S^{1} )$ such that 
\begin{equation*}
\limit{ k }{ \infty }{ \norm{ \left[ \psi_{ n_{0} , k  } \circ ( \phi_{ n_{0} , \infty  } \otimes \id_{ C( S^{1} ) } ) \right] ( x ) -  x  } } = 0
\end{equation*}
for all $x$ in $A_{ n_{0} } \otimes C( S^{1} )$.  

Set $\beta_{ n_{0} , k } = L_{ n_{0} , k } \circ \psi_{ n_{0}, k }$.  Suppose $a$ and $b$ are in $\mc{F}$.  Choose $x_{a}$ and $x_{b}$ in $\mc{F}_{ n_{0} }$ such that $\norm{ a - ( \phi_{ n_{0} , \infty } \otimes \id_{ C(S^{1}) } ) ( x_{a} ) } < \frac{\epsilon}{10}$ and $\norm{ b - ( \phi_{ n_{0} , \infty } \otimes \id_{ C(S^{1}) } ) ( x_{b} ) } < \frac{\epsilon}{10}$.  Then an easy computation shows that 
\begin{align*}
\norm{ \beta_{ n_{0} , k } ( a b ) - \beta_{ n_{0} , k } ( a  ) \beta_{ n_{0} , k } ( b ) } &< \frac{ \epsilon }{ 4 } + \norm{ \beta_{ n_{0} , k } ( x_{0} y_{0} ) - \beta_{ n_{0} , k } ( x_{0}  ) \beta_{ n_{0} , k } ( y_{0} ) },
\end{align*}  
where $x_{0} = ( \phi_{ n_{0} , \infty } \otimes \id_{ C(S^{1}) } ) ( x_{a} )$ and $y_{0} = ( \phi_{ n_{0} , \infty } \otimes \id_{ C(S^{1}) } )( x_{b} )$.  Since 
\begin{equation*}
\limit{ k }{ \infty }{ \norm{ \beta_{ n_{0} , k } ( x y ) - \beta_{ n_{0} , k } ( x ) \beta_{ n_{ 0 }, k } ( y ) } } = 0
\end{equation*}
for all $x$ and $y$ in $( \phi_{ n_{0} , \infty } \otimes \id_{ C(S^{1}) } ) ( A_{ n_{0} } \otimes C( S^{1} ) )$, we have that
\begin{equation*}
\limit{ k }{ \infty }{ \norm{ \beta_{ n_{0} , k } ( a b ) - \beta_{ n_{0} , k } ( a ) \beta_{ n_{0} , k } ( b ) } } \leq \frac{ \epsilon }{ 4 }
\end{equation*}
for all $a$ and $b$ in $\mc{F}$.  It is easy to check that  $\totalk ( \beta_{ n_{0} , k } \circ \iota ) \vert_{ \mc{P} } = \beta \vert_{ \mc{P} }$ and $\totalk ( \beta_{ n_{0} , k } )( [ 1_{A} \otimes z ] ) = \beta ( [ u ] )$ for $k$ sufficiently large.  

Choose $k_{0}$ in $\N$ such that 
\begin{equation*}
\norm{ \beta_{ n_{0} , k_{0} } (a b ) - \beta_{ n_{0} , k_{0} } ( a ) \beta_{ n_{0}, k_{0} } ( b ) } < \frac{ \epsilon }{ 2 }
\end{equation*}
for all $a$ and $b$ in $\mc{F}$.  Choose a finite subset $\mc{H}$ of $\proj ( B )$ such that for every $q$ in $\mc{Q}$, there exists $p_{q}$ in $\mc{H}$ such that $[ p_{q} ] =  \beta ( [ q ] )$.  Since $B$ is a unital, separable, amenable, simple tracally AF algebra, by Proposition 9.10 of \cite{linTR1}, there exists a sequence of contractive, completely positive, linear maps $\{ \gamma_{k} \}_{ k = 1 }^{ \infty }$ from $B$ to $A$ such that 
\begin{equation*}
\limit{ k } { \infty } { \norm{ \gamma_{k} ( x y ) - \gamma_{k}( x ) \gamma_{k} ( y ) } } = 0
\end{equation*}  
for all $x$ and $y$ in $B$ and 
\begin{equation*}
\totalk ( \gamma_{k} ) \vert_{ \mc{H} } = ( \beta^{-1} ) \vert_{ \mc{H} }
\end{equation*}
for all $k$.  Choose $k_{1}$ such that  
\begin{equation*}
\norm{ \gamma_{ k_{1} } ( x y ) - \gamma_{ k_{1} } ( x ) \gamma_{ k_{1} } ( y ) } < \frac{ \epsilon }{ 2 }
\end{equation*}
for all $x$ and $y$ in $\beta_{ n_{0} ,k_{0} } ( \mc{F} )$.

Set $\phi = \gamma_{ k_{1} } \circ \beta_{ n_{0} , k_{0} }$.  Then 
$\totalk ( \phi ) ( [ 1_{A} \otimes z ] ) = [ u ]$, $\totalk ( \phi \circ \iota ) \vert_{ \mc{P} } = \totalk( \id_{A} ) \vert_{ \mc{P} }$, and $\norm{ \phi ( x y ) - \phi ( x ) \phi ( y ) } < \epsilon$ for all $x$ and $y$ in $\mc{F}$.
\end{proof}

We next show that the contractive, completely positive, linear map obtained in the above lemma can be perturbed in such a way that the contractive, completely positive, linear map obtained by this small perturbation is ``close'' to $\id_{A}$ on the tracial state space of $A$. 

\begin{lemma}\label{lem:homtrace}
Let $A$ be a separable, unital, amenable, simple \cstar-algebra satisfying the Universal Coefficient Theorem.  Suppose $A$ is a tracially AI algebra.  Then for every $\epsilon$ in $\R_{ > 0 }$, finite subset $\mc{P}$ of $\proj ( A )$, finite subset $\mc{F}_{1}$ of $A_{\sa}$, and finite subset $\mc{F}_{2}$ of $A \otimes C( S^{1} )$, there exists a finitely generated subgroup $G$ of $K_{0} ( A )$ containing $[ 1_{A} ]$ such that the following holds:  for every unitary $u$ in $A$ with $[ u ] $ an element of $ H_{ [ 1_{A}] } ( G , K_{1} ( A ) )$, there exists a unital, completely positive, linear map $\psi$ from $A \otimes C( S^{1} )$ to $A$ such that 
\begin{enumerate}
\item $\psi$ is $\mc{F}_{2}$-$\epsilon$-multiplicative;
\item $\totalk ( \psi \circ \iota ) \vert_{ \mc{P} } = \totalk ( \id_{A} ) \vert_{ \mc{P} } $; 
\item $\totalk ( \psi ) ( [ 1_{A} \otimes z ] ) = [ u ]$; and
\item $\sup \setof{ | ( \tau \circ \psi \circ \iota ) ( a ) - \tau ( a ) | }{ \tau \in T( A ) } < \epsilon$ for all $a$ in $\mc{F}_{1}$.
\end{enumerate}  
\end{lemma}

\begin{proof}
Let $\delta$ and $\mc{G}$ be the quantities given in Lemma \ref{lem:pertunital} corresponding to $A \otimes C(S^{1})$, $\iota ( \mc{F}_{1} ) \cup \mc{F}_{2} \cup \{ 1_{A \otimes C(S^{1}) } \}$, $\iota( \mc{P} ) \cup \{ 1_{A} \otimes z  , 1_{ A \otimes C(S^{1}) } \}$, and $\frac{ \epsilon }{ 2 }$.  Let $\{ \mc{H}_{n} \}_{ n = 1 }^{ \infty }$ be an increasing sequence of finite subset of $A$ whose union is dense in $A$.  Note that we may assume that $\delta < \frac{ \epsilon }{ 2 }$.  Now, for each $n$, there exist a projection $p_{n}$ in $A$, a sub-\cstar-algebra $D_{n} = \bigoplus_{ i = 1 }^{ k(n) } \mrm{ M }_{ m( i, n) } ( C ( X_{ [ i, n ] } ) )$ of $A$, where $X_{ [ i, n ] } $ is either $\C$ or $[ 0 , 1 ]$ with $1_{ D_{n} } = p_{n} $, and a sequence of contractive, completely positive, linear maps $\{ L_{n} \}_{ n = 1 }^{ \infty }$ from $A$ to $D_{n}$ such that 
\begin{enumerate}
\item $\norm{ p_{n} x - x p_{n} } < \frac{ 1 }{ 2^{n} }$ for all $x$ in $\mc{H}_{n}$;
\item $\norm{ p_{n} x p_{n} - L_{n} ( x ) } < \frac{ 1 } { 2^{n} }$ for all $x$ in $\mc{H}_{n}$;
\item $\norm{ x - ( 1_{A} - p_{n} ) x ( 1_{A} - p_{n} ) - \psi_{n} ( x ) } < \frac{ 1 }{ 2^{n} }$ for all $x$ in $\mc{H}_{n}$ with $\norm{ x } \leq 1$; and
\item $\tau ( 1_{A} - p_{n} ) < \frac{1}{ 2^{n} }$ for all $\tau$ in $T(A)$.
\end{enumerate}

Note that $\limit{n}{ \infty }{ \norm{ L_{n} ( x y ) - L_{n} ( x ) L_{n} ( y ) } } = 0 $ for all $x$ and $y$ in $A$.  Denote the $i^{th}$ summand of $D_{n}$ by $D_{[n,i]}$.  Let $d_{ [ n,i]} = 1_{ D_{ [n,i] } }$.  Choose $n$ large enough such that $\frac{1}{2^{n}} < \frac{ \delta }{ 3}$.  Let $\mc{P}_{1}$ be a finite subset of $\proj ( A \otimes C(S^{1}) )$ such that $\mc{P}_{1}$ contains $\iota ( \mc{P} )$, $d_{ [ n , i ] } \otimes 1_{ C(S^{1}) }$, $p_{n} \otimes 1_{ C( S^{1} ) }$.  Choose a finite subset $\mc{X}_{2}$ of $A \otimes C(S^{1})$ such that $\mc{X}_{2}$ contains $\iota ( F_{1} ) \cup \mc{F}_{2}$ and the set
\begin{equation*}
\setof{ [ ( 1_{A} - p_{n} ) \otimes 1_{ C(S^{1}) } ] x [ ( 1_{A} - p_{n} ) \otimes 1_{ C(S^{1}) } ] }{ x \in \mc{F} }. 
\end{equation*}   

Let $G$ be the finitely generated subgroup of $K_{0} ( A )$ in Lemma \ref{lem:liftinghom} which corresponds to $\frac{1}{2^{n}}$, $\mc{P}$, and $\mc{X}_{2}$.  Suppose $[ u ]$ is an element of $H_{ [ 1_{A} ] } ( G , K_{1} ( A ) )$.  Then by Lemma \ref{lem:liftinghom}, there exists a contractive, completely positive, linear map $L$ from $A \otimes C(S^{1})$ such that 
\begin{enumerate}
\item $L$ is $\mc{X}_{2}$-$\frac{1}{2^{n}}$-multiplicative;
\item $\totalk ( L ) \vert_{ \mc{P}_{1} }$ is well-defined;
\item $\totalk ( L \circ \iota ) \vert_{ \mc{P} } = \totalk ( \id_{A} ) \vert_{ \mc{P} }$; and
\item $\totalk ( L ) ( [ 1_{A} \otimes z ] ) = [ u ]$ in $K_{1} ( A )$.
\end{enumerate}

Choose a projection $q_{n}$ in $A$ such that $ [ q_{n} ] = \totalk ( L ) ( \sum_{ i = 1}^{ k(n) } [ d_{ [ n, i ] } \otimes 1_{ C(S^{1} ) } ] )$.  Let $\mc{G}_{n}$ be a finite subset of $D_{n}$ such that $\mc{G}_{n}$ contains the generators of $D_{n}$.  Define $\gamma$ from $T ( A )$ to $T ( D_{n} )$ by $\gamma ( \tau ) = \frac{ 1 }{ \tau ( p_{n} ) } \tau \vert_{ D_{n} }$.  Using the same argument as in Proposition 9.7 of \cite{linTR1}, we get a $*$-homomorphism $h$ from $D_{n}$ to $q_{n} A q_{n}$ such that 
\begin{equation*}
\sup \setof{ | ( \tau \circ h ) ( g ) - \gamma ( \tau ) (g) ) | }{ \tau \in T(A) } < \frac{1}{ 2^{n} }
\end{equation*} 
for all $g$ in $\mc{G}_{n}$.  Define $\psi$ from $A \otimes C(S^{1})$ to $A$ by 
\begin{equation*}
\psi_{0} ( x ) = L ( [ (1_{A} - p_{n} ) \otimes 1_{ C( S^{1} ) } ] x [ (1_{A} - p_{n} ) \otimes 1_{ C( S^{1} ) } ] ) + ( h \circ L_{n} \circ \pi ) ( x ). 
\end{equation*}
Hence, $\psi_{0}$ is $\mc{X}_{2}$-$\delta$-multiplicative.  By construction, we have that 
\begin{enumerate}
\item $\totalk ( \psi_{0} ) ( [ 1_{A} \otimes z ] ) = \totalk ( L ) ( [ 1_{A} \otimes z ] ) = [ u ]$;
\item $\totalk ( \psi_{0} \circ \iota ) \vert_{ \mc{P} } = \totalk ( \id_{A} ) \vert_{ \mc{P} }$; and
\item for all $a$ in $\mc{F}_{1}$, 
\begin{equation*}
\sup \setof{ | ( \tau \circ \psi_{0} \circ \iota ) ( a ) - \tau ( a ) | }{ \tau \in T ( A ) } < \delta.
\end{equation*}
\end{enumerate}
By Lemma \ref{lem:pertunital}, there exists a contractive, completely positive, linear map $\beta$ from $A \otimes C(S^{1})$ to $A$ such that $\beta$ is $\mc{X}_{2}$-$\frac{ \epsilon } {2}$-multiplicative, $\beta ( 1_{A \otimes C( S^{1} ) } )$ is a projection in $A$, and 
\begin{equation*}
\norm{ \psi_{0} ( x ) - \beta ( x ) } < \frac{ \epsilon }{ 2 }
\end{equation*}
for all $x$ in $\mc{X}_{2}$.  Also, we have that $\totalk ( \beta \circ \iota ) \vert_{ \mc{P} } = \totalk ( \psi_{0} \circ \iota ) \vert_{ \mc{P} }$ and $\totalk ( \beta ) ( [ 1_{A} \otimes z ] ) = \totalk ( \psi_{0} ) ( [ 1_{A} \otimes z ] )$.  Hence, $\totalk ( \beta ) ( [ 1_{ A \otimes C(S^{1}) }  ] )  = [ 1_{A} ]$.  Since $A$ has stable rank one, there exists a unitary $w$ in $A$ such that $w \beta ( 1_{ A \otimes C(S^{1}) } ) w^{*} = 1_{A}$. 

Set $\psi = \innerauto ( w ) \circ \beta$.  It is easy to check that $\psi$ is the desired unital, completely positive, linear map from $A \otimes C(S^{1})$ to $A$.  
\end{proof}

Let $A$ be a unital \cstar-algebra. Suppose that $p$ is a projection in $A$, $a$ is an element of $A$ with $ \norm{ a}  \leq 1$, and
\begin{equation*}
\norm{ a^*a-p } < \frac{1}{16} \ \mrm{and} \ \norm{ aa^*-p } < \frac{1}{16}.
\end{equation*}
A standard computation shows
that
\begin{equation*}
\norm{ pap - ap } < \frac{ 3 }{ 16 } \ \mrm{and}\  \norm{ pa - pap } < \frac{ 3 }{ 16 }.
\end{equation*}
Also
$
\norm{ pa - a } < \frac{1}{2}.
$
Set $b=pap.$ Then
\begin{equation*}
\norm{ b^*b - p } \leq \norm{ pa^*ap-pa^*a } + \norm{ pa^*a-p } < \frac{ 1 }{ 16 } + \frac{ 1 }{ 16 } = \frac{1}{8}.
\end{equation*}
So
\begin{equation*}
\norm{ (b^*b)^{-1} - p } < \frac{ \frac{1}{8} }{ 1 - \frac{1}{8} } = \frac{1}{7} \ \mrm{and} \ \norm{ |b|^{-1} - p } < \frac{2}{7},
\end{equation*}
where the inverse is taken in $pAp.$
Set $v=b|b|^{-1}.$
Then $v^*v=p=vv^*$ and
\begin{equation*}
\norm{ v - b } < \frac{2}{7}.
\end{equation*}
We denote
$v$ by $\tilde{a}$.  Note that if $v_{1}$ is another unitary in $pAp$ with
$\norm{ v_{1} - b } < \frac{1}{3},$ then $[ v_{1} ] = [ v ]$ in $U(pAp)/U_0(pAp).$

Suppose that $L$ from $A$ to $B$ is a $\mc{G}$-$\delta$-multiplicative, contractive, completely positive, linear map, $u$ is a normal partial isometry, and a projection
$p$ in $B$ is given so that
\begin{equation*}
\norm{ L(u^*u)-p } < \frac{1}{32}.
\end{equation*}
We define  $\tilde{L}$ as follows.
Let $L(u)=a.$ With $\delta$ chosen to be sufficiently small and $\mc{G}$ chosen to be sufficiently large,
we denote by $\tilde{ L }(u)$ the normal partial isometry
(unitary in a corner)
$v = \tilde{a}$ defined above. This notation will be used later.
Note also, if $u$ in $U_{0} (A)$, then with sufficiently large
$\mc{G}$ and sufficiently small $\delta$, we may assume
that $\tilde{L}(u)$ is an element of $U_{0}(B)$.

Let $\phi$ from $ [ 0 , a ]$ to $X$ be a continuous map, where $X$ is a normed space.  Let $\mc{P} = \{ 0 = t_{0} < t_{1} < \cdots < t_{n} = a \}$ be a partition of $[ 0 , a ]$.  Set $L ( \phi ) ( \mc{P} ) = \sum_{ i = 1}^{n} \norm{ \phi ( t_{i} ) - \phi ( t_{i-1} ) }$.  The \emph{length} of $\phi$, denoted as $L (\phi)$, is defined to be the supremum of $L ( \phi ) ( \mc{P} )$ over all partitions $\mc{P}$ of $[0,a]$.  

If $A$ is a unital \cstar-algebra and $u$ is an element of $U(A)_{0}$, then $\mrm{cel}(u)$ is defined to be 
\begin{equation*}
\mrm{cel} ( u ) = \inf \setof{ L(\phi) }{ \ftn{ \phi }{ [ 0 ,1 ] }{ U(A)_{0} } , \ \phi ( 0 ) = 1_{A} , \ \phi ( 1 ) = u }. 
\end{equation*}

\begin{definition}\label{def:commutatorsubgrp}
Let $A$ be a unital \cstar-algebra.  Let $CU(A)$ be the closure of the commutator subgroup of $U(A)$.  Clearly, that the commutator subgroup forms a normal subgroup of $U(A)$.  Also note that $U(A) / CU ( A )$ is commutative.  If the natural map from $U(A) / U_{0} ( A )$ to $K_{1} ( A )$ is injective, then $CU ( A )$ is a normal subgroup of $U_{0} ( A )$.  If $u$ is an element of $U(A)$, we will denote the image of $u$ in $U ( A ) / CU ( A )$ by $\bar{u}$, and if $F$ is a subgroup of $U(A)$, then $\bar{F}$ will denote the image of $F$ in $U(A) / CU ( A )$.  

If $\bar{u}$ and $\bar{ v}$ are elements of $U(A) / CU(A)$, define
\begin{equation*}
\mrm{dist} ( \bar{u}, \bar{v} ) = \inf \setof{ \norm{ x - y  } }{  x, y \in U(A)\,\,\, {\rm such\,\,\, that}\,\,\, \bar{x} = \bar{u}, \, \bar{y} = \bar{v} }.
\end{equation*}
If $u$ and $v$ are elements of $U(A)$, then
$\mrm{dist} ( \bar{u}, \bar{v} ) = \inf \setof{ \norm{ uv^*- x } } { x \in CU(A) }$.  Let $g = \prod_{ i = 1}^{n} a_{i} b_{i} a_{i}^{-1} b_{i}^{-1}$, where $a_{i}$ and $b_i$ are elements of $U(A)$.  Let $\mc{G}$ be a finite subset of $A$, $\delta$ in $\R_{ > 0}$, and $L$ from $A$ to $B$ be a $\mc{G}$-$\delta$-multiplicative, contractive, completely positive, linear map, where $B$ is a
unital \cstar-algebra. From the paragraphs before the definition, for $\epsilon$ in $\R_{ > 0}$ if $\mc{G}$ is sufficiently large
and $\delta$ is sufficiently small,
\begin{equation*}
\norm{ L( g ) - \prod_{ i = 1}^{n} x_{i} y_{i} (x_{i})^{-1}( y_{i})^{-1} } < \frac{ \epsilon }{2},
\end{equation*}
where  $x_{i}$ and $y_{i}$ are in $U(B)$.  Thus, for any $g$ in $CU(A)$, with sufficiently large $\mc{G}$ and
sufficiently small $\delta,$
\begin{equation*}
\norm{ L(g) - u } < \epsilon
\end{equation*}
for some $u$ in $CU(B)$.  Moreover, for any finite subset
$\mc{U}$ of $U(B)$ and subgroup
$ F$ of $U(B)$ generated by $\mc{U}$,
and $\epsilon$ in $\R_{ > 0}$, there exist a finite subset $\mc{G}$ and
$\delta$ in $\R_{ > 0}$ such that, for any $\mc{G}$-$\delta$-multiplicative, contractive, completely positive, linear map $L$ from $A$ to $B$, we have that $L$ induces a homomorphism, $L^{\ddag}$ from $\bar{F}$ to $U(B) / CU(B)$ such that $\mrm{dist} ( \overline{ \tilde{L} (u) } , L^{\ddag} ( \bar{u} ) ) < \epsilon$
for all $u$ in $\mc{U}$.  

If $\phi$ from $A$ to $B$ is  a $*$-homomorphism, then $\phi$ induces a continuous homomorphism 
$\phi^{\ddag}$ from $U(A) / CU(A)$ to $U(B) / CU(B)$.
\end{definition} 

We are now ready to prove Theorem \ref{thm:apuniteq} in the tracially AI case.  The proof follows the same line of argument as in the proof of Theorem 10.4 of \cite{linTR1}.

\noindent\emph{Proof of Theorem \ref{thm:apuniteq}: The tracially AI case} Let $\epsilon$ be an element of $\R_{ > 0}$ and $\mc{F}$ be a finite subset of $A$.  Suppose $A$ is a finite dimensional simple tracially AI algebra.  Then $A = \msf{M}_{n}$.  Let $G = K_{0} ( A ) = \Z$.  Since $K_{1} ( A ) = 0$, we can take $\Phi$ to be the canonical projection of $A \otimes C(S^{1})$ to $A$ and $w = 1_{A}$.  Then $K_{1} ( \Phi )( [ 1_{A} \otimes z ] ) = 0 = [ u ]$ in $K_{1} ( A )$ for all unitaries $u$ in $A$ and 
$\Phi \circ \iota = \id_{A}$.

Suppose $A$ is an infinite dimensional \cstar-algebra.  By Lemma 10.9 and Theorem 10.10 of \cite{linTR1}, $A$ is isomorphic to $\dirlim ( A_{n} , \phi_{n, n+1} )$, with $A_{n} = \bigoplus_{ j = 1}^{ m(n) } A_{ [ n , j ] }$ is as described in Theorem 10.1 of \cite{linTR1} such that $\phi_{ n , m}$, $\phi_{ n , \infty }$ and $K_{*} ( \phi_{ n , m } )$ are injective maps for all $m > n$.  For notational convenience, we may identify $A_{n}$ with $\phi_{ n , n+1} ( A_{n} )$.  We will use this identification without further warning.  

Define $\mbf{ L }$ from $U(A)$ to $\R_{ \geq 0}$ as follows:  if $u$ is an element of $U_{0} ( A )$, then $\mbf{ L } ( u ) = 2 \mrm{ cel } ( u ) + 8 \pi + \frac{ \pi }{ 16 }$; if $u$ is an element of $U( B ) - U_{0} ( B )$ and if there exists a positive integer such that $u^{k}$ is an element of $U_{0} ( A )$, then $\mbf{L} ( u ) = 16 \pi + \frac{ \mrm{ cel } ( u^{ k(u) } ) }{ 16 } + \frac{ \pi }{ 16 }$ where $k(u)$ is the order of $[ u ]$ in $U ( A ) / U_{0} ( A )$; and if $u$ is an element of $U(A) - U_{0} ( A )$ and $[ u ]$ does not have finite order in $U(A) / U_{0} ( A )$, then $\mbf{L} ( u ) = 16 \pi + \frac{ \pi }{ 16 }$.  

Note that we may assume that $\epsilon$ is less than $\frac{ \pi }{ 128 }$.  Given $A$, $\mbf{L}$, $\epsilon$, and $\mc{F}$, Theorem 8.6 of \cite{linTR1} provides us with an element $\delta_{1}$ of $\R_{ > 0}$, an element $n$ of \N, a subset finite $\mc{P}$ of $\proj(A)$, and a finite subset $\mc{S}$ of $A$.  Theorem 8.6 of \cite{linTR1} also provides us with mutually orthogonal projections $q, p_{1}, \dots, p_{n}$ such that for each $i$, we have that $q$ is Murray-von Neumann equivalent to a subprojection of $p_{i}$ and $p_{i}$ is Murray-von Neumann equivalent to $p_{1}$ and there exist a sub-\cstar-algebra $C_{1}$ in $\mc{I}$ with $1_{C_{1}} = p_{1}$ and unital, contractive, completely positive, linear maps $h_{0}$ from $A$ to $q A q$ and $h_{1}$ from $A$ to $C_{1}$ such that 
\begin{enumerate}
\item $h_{0}$ and $h_{1}$ are $\mc{S}$-$\frac{ \delta_{1} }{ 4 }$-multiplicative;
\item $h_{0}(x) = q x q$; and
\item $\norm{ x - ( h_{0} ( x ) \oplus \underbrace{h_{1} ( x ) \oplus \cdots \oplus h_{1} ( x )}_{n} ) } < \frac{ \delta_{1} }{ 16 }$ for all $x$ in $\mc{S}$.  
\end{enumerate}
Set $C = M_{ n } ( C_{1} )$ which is identified as a sub-\cstar-algebra of $(1_{A} - q ) A ( 1_{A} - q )$.  With this choice of $C_{1}$, Theorem 8.6 of \cite{linTR1} provides us with a finite subset $\mc{G}_{0}$ of $A$, a finite subset $\mc{P}_{0}$ of projections in $\msf{M}_{ \N } ( C )$, a finite subset $\mc{H}$ of $A_{\sa}$, elements $\delta_{0}$ and $\sigma$ of $\R_{ > 0}$.  Set $\delta$ equal to the minimum of $\delta_{0}$ and $\delta_{1}$.  

Note that we may assume that $\mc{P}_{0}$ contains $1_{A} - q$ and contains at least one minimal projection of each summand of $C$.  We also may assume that for each $u$ in $U( A ) \cap \mc{P}$, $u$ has the form $q u q \oplus ( 1_{A} - q ) u ( 1_{A} - q )$ where $q u q$ is a unitary in $q A q$ and $(1_{A} - q) u (1_{A} - q )$ is a unitary in $C$.  Note that there exists $n_{0}$ in $\N$ and a projection $q_{0}$ in $A_{n_{0}}$ such that $\phi_{ n_{0} , \infty }( q_{0} )$ is unitarily equivalent to $q$ via a unitary $v$ in the connected component of the identity.  Hence, by conjugating $\phi_{n, \infty }$ by $\innerauto(v)$, we may assume that $\phi_{ n_{0} , \infty } ( q_{0} ) = q$.  Also, we may assume that $q u q$ is in $\phi_{ n_{0} , \infty } ( A_{n_{0}} )$.  

Set $\mc{U} = \setof{ q u q }{ u \in U(A) \cap \mc{P} }$.  The subgroup of $ U( q A q )$ generated by $\mc{U}$ will be denoted by $F$ and $\overline{F}$ will denote the image of $F$ in $U ( q A q ) / CU ( q A q )$.  By Theorem 6.6(3) of \cite{linTR1}, $\overline{F} = \overline{F} \cap \left( U_{0} ( q A q ) / CU ( q A q ) \right) \oplus \overline{F}_{0} \oplus \overline{F}_{1}$, where $\overline{F}_{0}$ is a torsion group and $\overline{F}_{1}$ is a free group.  If $\kappa$ denotes the homomorphism from $U( q A q ) / CU( q A q)$ to $K_{1} ( q A q )$, then $\kappa ( \overline{ F }_{1} )$ is isomorphic to $\overline{F}_{1}$.  

Suppose $\mc{U}$ has the following decomposition:
\begin{equation}\label{decomp}
\mc{U} = \mc{U}_{0} \cup \mc{U}_{1}
\end{equation}
such that $\overline{ \mc{U}_{0} }$ generates $\overline{F} \cap \left( U_{0} ( q A q ) / CU ( q A q ) \right) \oplus \overline{F}_{0}$ and $\overline{ \mc{U}_{1} }$ generates $\mc{F}_{1}$.  It turns out that we can reduce the general case to this case since by Lemma 6.9 of \cite{linTR1}, modulo unitaries in $CU( q A q )$, this decomposition can be made with the cost of no more than $8\pi$ in the estimation of the exponential length.  Also, choosing a larger $n_{0}$ if necessary, we may assume that $\mc{U}_{0}$ and $\mc{U}_{1}$ are subsets of $q \phi_{ n_{0} , \infty } ( A_{n_{0}} ) q$.

For the quantities $\iota ( \mc{S} \cup \mc{G}_{0} \cup \mc{F} ) \cup \{ 1_{A} \otimes z , 1_{A} \otimes z^{*} \}$, $\mc{H}$, $\mc{P} \cup \mc{P}_{0}$, and $\delta_{2} = \frac{ \min\{ \delta , \sigma, \epsilon \} }{ 100 }$, Lemma \ref{lem:homtrace} provides us with a finitely generated subgroup $G$ of $K_{0} ( A )$ containing $[1_{A}]$.  Suppose $u$ is a unitary in $A$ such that $[ u ]$ is an element of $H_{ [ 1_{A} ] } ( G , K_{1} ( A ) )$.  By Lemma \ref{lem:homtrace}, there exists a unital, completely positive, linear map $\psi$ from $A \otimes C(S^{1})$ to $A$ such that $\psi$ is $\iota( \mc{S} \cup \mc{G}_{0} \cup \mc{F} ) \cup \{ 1_{A} \otimes z , 1_{A} \otimes z^{*} \}$-$\delta_{2}$-multiplicative and 
\begin{enumerate}
\item $\totalk ( \psi \circ \iota ) \vert_{ \mc{P} \cup \mc{P}_{0} } = \totalk ( \id_{A} ) \vert_{ \mc{P} \cup \mc{P}_{0} }$;
\item $\sup_{ \tau \in T(A) } \setof{ | ( \tau \circ \psi \circ \iota )( x ) - \tau ( x ) | }{ \tau \in T(A) } < \delta_{2}$ for all $x$ in $\mc{H}$; and 
\item $\totalk ( \psi ) ( [ 1_{A} \otimes z ] ) = [ u ]$ in $K_{1} ( A )$.
\end{enumerate}
Since $A_{n}$ and $A$ are separable, amenable \cstar-algebras and $A = \dirlim ( A_{n} , \phi_{ n , n+1} )$, there exists a sequence of unital, completely positive, linear maps $\{ \ftn{ \mu_{ n } }{ A }{ A_{n} } \}_{ n = 1}^{ \infty }$ such that $\limit{ n }{ \infty }{ \norm{ \mu_{n} ( x y ) - \mu_{n} ( x ) \mu_{n} ( y ) } } = 0$ for all $x$ and $y$ in $A$ and $\norm{ ( \phi_{ n , \infty } \circ \mu_{n} ) ( a ) - a } = 0$ for all $a$ in $A$.  Therefore we may choose $n_{1} \geq n_{0}$ such that 
\begin{enumerate}
\item $\mu_{ n_{1} } \circ \psi$ is $\iota( \mc{S} \cup \mc{G}_{0} \cup \mc{F} ) \cup \{ 1_{A} \otimes z , 1_{A} \otimes z^{*} \}$-$2 \delta_{2}$-multiplicative;
\item $\norm{ \phi_{ n_{1} , \infty } \circ \mu_{n_{1}} \circ \psi ( x ) - \psi( x ) } < 2 \delta_{2}$ for all $x$ in $\iota( \mc{S} \cup \mc{G}_{0} \cup \mc{F} ) \cup \{ 1_{A} \otimes z , 1_{A} \otimes z^{*} \}$;
\item $\totalk ( \phi_{ n_{1} , \infty } \circ \mu_{ n_{1} } \circ \psi \circ \iota ) \vert_{ \mc{P} \cup \mc{P}_{0} } = \totalk ( \psi \circ \iota ) \vert_{ \mc{P} \cup \mc{P}_{0} }$; 
\item $\totalk ( \phi_{ n_{1} , \infty } \circ \mu_{ n_{1} } \circ \psi )( [ 1_{A} \otimes z ] ) = \totalk ( \psi )( [ 1_{A} \otimes z ] ) = [ u ]$ in $K_{1} ( A )$; and 
\item $\sup \setof{ | ( \tau \circ \phi_{ n _{1} , \infty } \circ \mu_{ n_{1} } \circ \psi \circ \iota )( a ) - \tau ( a ) | } { \tau \in T(A) } < 2 \delta_{2}$ for all $a$ in $\mc{H}$.
\end{enumerate}  

Let $B_{1} = q A_{ n_{1} } q$.  Since $A$ is simple, it is known and easy to see
that, by choosing possibly a large $n_{1}$, we may assume that the
rank of $q$ at each point is at least $6$ (in $A_{ n_{1} }$).  Note that we have assumed that
$q$ is an element of $A_{ n_{1} }$.  So $B_{1}$ is a corner of $A_{ n_{1} }$.  By construction and the fact that each $K_{*} ( \phi_{ n , \infty } )$ is injective,
we have that  $\totalk ( \mu_{ n_{1} } \circ \psi \circ \iota ) ( [ q ] ) = [ q ] $.  By conjugating $\mu_{ n_{1} } \circ \psi \circ \iota $ by some unitary
$w$ if necessary, we may assume that
\begin{equation*}
\norm{ ( \mu_{ n_{1 } } \circ \psi \circ \iota ) (q) - q } < \frac{ \delta }{ 4 }.
\end{equation*}
Define $\Lambda( b ) = a q [  ( \mu_{ n_{1 } } \circ \psi \circ \iota )(qbq) ] q a$
(where $a= [ q ( \mu_{ n_{1 } } \circ \psi \circ \iota )(q) q ]^{-1/2}$) for $b$ in $q A q$.
Note that
\begin{equation*}
\norm{ \Lambda - ( \mu_{ n_{1 } } \circ \psi \circ \iota ) \vert_{ q A q } } < \frac{ \delta }{ 2}.
\end{equation*}

Write $A_{ n_{1} } = \oplus_{ k = 1}^{m} A_{ [ n(k) , n_{1} ] }$, where each $A_{ [ n(k) , n_{1} ] } $
has the form $C^{(k)}$ as described in Definition 7.1 of \cite{linTR1}.
According to this direct sum decomposition, we may write
$q = q_{1} \oplus q_{2} \oplus \cdots \oplus q_{l}$ with $0 \leq l \leq m$
and $q_{k} \neq 0$, for $1 \leq k \leq l$.  Choose an integer $N_{1} > 0$ such that $N_{1} [ q_{k} ] \geq 3 [ 1_{ A_{ [ n(k) , n_{1} ] } } ]$ for $k \leq l$.  Note that we may assume that $q_{ k }$ has rank at least $6$.  By applying an inner automorphism, we may assume that
$\oplus_{ k = 1}^{l} A_{ [ n(k), n_{1} ] }$ is a hereditary sub-\cstar-algebra of $\msf{M}_{ N_{1} } ( B_{1} )$.  Since $F_{1}$ is finitely generated, with sufficiently large $n_{1}$,
we obtain (see Definition 6.2 of \cite{linTR1}) a homomorphism $j$ from $\bar{F}_{1}$ to $U( B_{1} ) / CU( B_{1} )$
such that $\phi_{ n_{1} }^{\ddag} \circ j = \rm{id}_{ \bar{F}_{1} }$, where $\phi_{ n_{1} } = \phi_{ n_{1} , \infty }$.  Then (since the canonical map from $K_{1} ( A_{ n_{1} } )$  to $K_{1} ( A )$ is injective),
\begin{equation*}
\kappa_{1} \circ \phi_{ n_{1} }^{\ddag} \circ ( \mu_{ n_{1 } } \circ \psi \circ \iota )^{\ddag} \vert_{ \overline{ F }_{1}  } = \kappa_{1} \circ (\phi_{ n_{1} } )^{\ddag} \circ j = (\kappa_{1} )_{ \overline{F}_{1} },
\end{equation*}
where $\kappa_{1}$ from $ U( q A q ) / CU( q A q )$ to $K_{1}( q A q )$ is the quotient map.
Note that $K_{1}( q A q ) = K_{1}( A )$.  Let $\Delta$ be $\delta( \frac{\epsilon}{4} )$ as
described in Lemma 7.5 of \cite{linTR1}. We may assume that $\Delta < \frac{ \sigma }{ 4}$.  To
simplify notation, without loss of generality, we may assume that
$\phi_{ n_{1} } ( q ) = q$.  By the assumption on $A$, we may write that
$\phi_{ n_{1} } \vert_{ B_{1} } = ( \phi_{ n_{1} } )_{0} \oplus ( \phi_{ n_{1} } )_{1}$, where
\begin{enumerate}
\item $\tau( ( \phi_{ n_{1} } )_{0} ( 1_{ B_{1} } ) ) < \frac{ \Delta }{ 2(N_1+1)^2 }$ for all $\tau$ in $T(A)$ and

\item $( \phi_{ n_{1} } )_{0}$ is homotopically trivial (but nonzero).
\end{enumerate}
(see Theorem 10.1 of \cite{linTR1})

It follows from Lemma 7.5 of \cite{linTR1} that there exists a $*$-homomorphism $h$ from $B_{1}$ to $e_{0} A e_{0}$ such that
\begin{itemize}
\item[(i)] $\totalk( h ) =  \totalk ( ( \phi_{ n_{1} } )_{0} )$ in $\Hom_{ \Lambda } ( \totalk(B_{1}), \totalk(A) )$ and

\item[(ii)] $( \phi_{ n_{1} }^{\ddag} \circ j ( \bar{w} ) )^{-1} ( h \oplus ( \phi_{ n_{1} } )_1)^{\ddag}( \Lambda^{\ddag} ( \bar{w} ) )
= \overline{ g_{w} }$, where $g_{w}$ is an element of $U_0( q A q )$ and $\mrm{cel} ( g_{w} ) < \frac{ \epsilon }{ 4}$ (in
$U( q A q)$) for all $w$ in $\mc{U}_1.$
\end{itemize}

Recall that we have assumed that $A_{ n_{1} }$ is a sub-\cstar-algebra of $\msf{M}_{ N_{1} } ( B_{1} )$.  Set
\begin{equation*}
t = ( ( h \oplus ( \phi_{ n_{1} } )_{1} ) \otimes \rm{id}_{ \msf{M}_{ N_{1} } } ) \vert_{ \oplus_{ k = 1 }^{l} A_{ [ n(j) ,n_{1} ] } }.
\end{equation*}
and $\Psi = t \oplus ( \phi_{ n_{1} } ) \vert_{ \oplus_{ k = l + 1 }^{m}  A_{ [ n(j), n_{1} ] } }$.
Let $\Phi = \Psi \circ \mu_{ n_{1} } \circ \psi$.  It is clear that (since $\Delta < \frac{ \sigma }{ 4}$)
\begin{enumerate}
\item $\totalk ( \Phi \circ \iota ) \vert_{ \mc{P} \cup \mc{P}_{0} }= \totalk ( \phi_{ n_{1} , \infty } \circ \mu_{n_{1}} \circ \psi \circ \iota ) \vert_{ \mc{P} \cup \mc{P}_{0} } = \totalk ( \psi \circ \iota) \vert_{ \mc{P} \cup \mc{P}_{0} }$;
\item $\totalk ( \Phi \circ \iota )( [ 1_{A} \otimes z ] ) = \totalk ( \phi_{ n_{1} ,\infty } \circ \mu_{ n_{1} } \circ \psi )( [ 1_{A} \otimes z ] ) = \totalk ( \psi )( [ 1_{A} \otimes z ] ) = [ u ]$ in $K_{1} ( A )$; and 
\item $\sup \setof{ | ( \tau \circ \Phi \circ \iota ) ( a ) - ( \tau \circ \phi_{ n_{1} , \infty } \circ \mu_{ n_{1} } \circ \psi )( a ) | }{ \tau \in T(A) } < \frac{ \sigma }{ 2 }$ for all $a$ in $A$.
\item for all $w$ in $\mc{U}_{1}$, by (ii) we have that
\begin{equation*}
\mrm{cel} ( w^{*} \widetilde{(\Phi \circ \iota) } ( w ) ) < 8 \pi + \frac{ \epsilon }{ 2 } \ \mrm{in} \ U ( q A q ); \ \mrm{and}
\end{equation*}
\item for $w$ in $\mc{U}_{0}$, by Lemma 6.8, Theorem 6.10, and Lemma 6.9 of \cite{linTR1} we have that
\begin{equation*}
\mrm{cel}( w^{*} \widetilde{ (\Phi \circ \iota) } ( w ) ) < 
\begin{cases}
2 \mrm{cel} ( w ) + \frac{ \pi }{ 64 }, &\text{in $U(q A q)$ if $[ w ] = 0$ in $K_{1} ( A )$}\\
8 \pi + \frac{ 2 \mrm{cel}( w^{ k(w) } ) }{ k(w) } + \frac{ \pi }{ 16 }, &\text{in $U(q A q)$ if the order of $[w]$ is $k(w)$}
\end{cases}
\end{equation*}
\end{enumerate}
Therefore, even after we add $8 \pi$ for the decomposition of $\overline{F}$ as in (\ref{decomp}), we get 
\begin{equation*}
\mrm{cel}( \id_{A} ( h_{0} ( u )^{*} ) \widetilde{ (\Phi \circ \iota) }( h_{0} ( u ) ) ) < \mbf{L}(u) \ \mrm{in} \ U( q A q )
\end{equation*}
for all $u$ in $U( A ) \cap \mc{P}$.  Since we also have that
\begin{equation*}
\totalk ( \Phi \circ \iota ) \vert_{ \mc{P} \cup \mc{P}_{0} } = \totalk ( \id_{A} ) \vert_{ \mc{P} \cup \mc{P}_{0} }
\end{equation*}
and
\begin{equation*}
\sup\setof{ | ( \tau \circ \Phi \circ \iota )( a ) - \tau ( a ) | }{ \tau \in T(A) } < \sigma
\end{equation*}
for all $a$ in $\mc{H}$,  we can apply Theorem 8.6 of \cite{linTR1} to $\Phi \circ \iota$ and $\id_{A}$ to get a unitary $w$ in $A$ such that 
\begin{equation*}
\norm{ w (\Phi \circ \iota )( a ) w^{*} - a } < \epsilon 
\end{equation*}
for all $a$ in $\mc{F}$.  Note that $\Phi$ is a unital, contractive, completely positive, linear map such that $\Phi$ is $\iota( \mc{F} ) \cup \{ 1_{A} \otimes z , 1_{A} \otimes z^{*} \}$-$\epsilon$-multiplicative and $\totalk ( \Phi )( [ 1_{A} \otimes z ] ) = [ u ]$ in $K_{1} (A)$. 

\begin{remark}
If, in the above theorem, we assumed that $A$ is tracially AF or $A$ is tracially AI with $K_{1} (A)$ a torsion group, then the proof of the above theorem would be much easier.  In these two cases, total $K$-theory and traces are enough to determine when two almost multiplicative linear map are approximately unitarily equivalent (see Theorem 6.3.3 of \cite{booklin} and Theorem 8.7 of \cite{linTR1}).  Hence, one does not need to control the exponential length, so the proof is much easier and shorter.
\end{remark}

\section{The automorphism group a simple \cstar-algebra}\label{auto}

For any \cstar-algebra $A$, we will denote the subgroup of $\Hom_{ \Lambda } ( \totalk ( A ) , \totalk ( A ) )$ consisting of all elements $\alpha$ such that $\alpha$ is an isomorphism which sends $[ 1_{A} ]$ to $[ 1_{A} ]$ and $\alpha \vert_{ K_{0} ( A ) }$ and $\alpha^{-1} \vert_{ K_{0} ( A ) }$ are positive homomorphisms by $\mrm{ Aut }_{ \Lambda } ( \totalk ( A ) )_{ + , 1}$. 

\begin{theorem}\label{thm:extgrps}
Let $A$ be a simple, unital, separable, amenable \cstar-algebra satisfying the Universal Coefficient Theorem.  If $A$ is a purely infinite \cstar-algebra or a tracially AF algebra, then
\begin{equation*}
\xymatrix{
0 \to \overline{ \mrm{Inn} } ( A ) \to \mrm{Aut} ( A ) \to \mrm{ Aut }_{ \Lambda } ( \totalk ( A ) )_{ + , 1 } \to 0
}
\end{equation*}
is exact and $\overline{ \mrm{Inn} } ( A )$ is an extension of a totally disconnected abelian topological group by a simple topological group.

Moreover, if there exists a sequence $\{ G_{n} \}_{ n = 1}^{ \infty }$ of finitely generated subgroup of $K_{0} ( A )$ whose union is $K_{0} ( A )$ and each $G_{n}$ contains $[ 1_{A} ]$ such that $H_{ [ 1_{A} ] } ( G_{n} , K_{1} ( A ) ) = K_{1} ( A )$ for all $n$ sufficiently large, then $\overline{ \mrm{ Inn } } ( A ) = \overline{ \mrm{Inn} }_{0} ( A )$.  Consequently, $\overline{ \mrm{ Inn } } ( A )$ is a simple topological group.  
\end{theorem}

\begin{proof}
By Theorem 6.7 of \cite{sepBDF} for the purely infinite case and by Theorem 6.3.3 of \cite{booklin} for the tracially AF case, $0 \to \overline{ \mrm{Inn} } ( A ) \to \mrm{Aut} ( A ) \to \mrm{ Aut }_{ \Lambda } ( \totalk ( A ) )_{ + , 1 } $ is exact.  We will now show that every element of $\mrm{ Aut }_{ \Lambda } ( \totalk ( A ) )_{ + ,1 }$ lifts to an automorphism.    Suppose $\alpha$ is an element of $\mrm{ Aut }_{ \Lambda } ( \totalk ( A ) )_{ + ,1}$.  By Theorem 1.1 of \cite{morphTAF} (in the tracially AF case) and by Theorem 6.7 of \cite{sepBDF} (in the purely infinite case), there exist unital $*$-homomorphisms $\gamma$ and $\beta$ from $A$ to $A$ such that $\gamma$ induces $\alpha$ and $\beta$ induces $\alpha^{-1}$. 

Let $\{ \mc{F}_{n} \}_{ n = 1}^{ \infty }$ be an increasing sequence of finite subsets of $A$ whose union is dense in $A$.  By Theorem 6.3.3 of \cite{booklin} (for the tracially AF case) and by Theorem 6.7 of \cite{sepBDF} (for the purely infinite simple case), there exists a sequence of unitaries $\{ w_{n} \}_{ n = 1}^{ \infty }$ in $A$ such that 
\begin{equation*}
\norm{ ( \beta_{n + 1} \circ \gamma_{n} )( x ) - x } < \frac{1}{ 2^{n} } \ \mrm{and} \ \norm{ ( \gamma_{n+1 } \circ \beta_{ n } )( x ) - x } < \frac{ 1 }{ 2^{n+1} }
\end{equation*}
for all $x$ in $\mc{F}_{ n}$, where $\gamma_{n} = \innerauto ( w_{2n}) \circ \gamma$ with $\gamma_{0} = \gamma$ and $\beta_{n} = \innerauto ( w_{ 2n +1} ) \circ \beta$.  So we have a two-sided approximately intertwining diagram
\begin{equation*}
\xymatrix{
A \ar[rd]_{ \gamma_{1} } \ar[rr]^{ \id_{A }}		&  & A \ar[rd]^{  \gamma_{2} }	\ar[rr]^{ \id_{A} }	& & A \ar[rr]  & & & \cdots & A					\\
	& A \ar[ru]_{ \beta_{1} } \ar[rr]_{ \id_{A} } & & A \ar[rr]_{ \id_{A} } \ar[ur]_{ \beta_{2} } & &  & & \cdots & A
	}
\end{equation*}
Therefore, we have an automorphism $\phi$ from $A$ to $A$ which induces $\alpha$.  

By Corollary 2.5 of \cite{EllRorAuto}, $\overline{ \mrm{Inn} }_{0}( A )$ is a simple topological group and by Theorem \ref{thm:inniso} and Theorem \ref{thm:propC}, $\frac{ \overline{ \mrm{Inn} } ( A ) }{ \overline{ \mrm{Inn} }_{0}( A ) }$ is a totally disconnected abelian topological group.  The last statement of the theorem follows from Theorem \ref{thm:inniso} and Theorem \ref{thm:propC}.
\end{proof}

\begin{theorem}\label{thm:minimalDS}
Let $X$ be an infinite compact metric space with finite covering dimension and let $\alpha$ be a minimal homeomorphism of $X$.  Denote the transformation group \cstar-algebra associated to $( X, \alpha )$ by $C^{*}( \Z , X, \alpha)$.  Suppose that the image $K_{0} ( C^{*} ( \Z , X , \alpha ) )$ under the map $\rho_{ C^{*} ( \Z , X , \alpha ) }$ is dense in $\mrm{Aff} ( T ( C^{*} ( \Z , X , \alpha ) ) )$.  Then 
\begin{enumerate}
\item $C^{*} ( \Z , X , \alpha )$ is a unital, separable, amenable, simple, tracially AF algebra which satisfies the Universal Coefficient Theorem;

\item $\overline{ \mrm{Inn} }_{0} (  C^{*}( \Z , X, \alpha) )$ is a simple topological group;

\item $\frac{ \overline{ \mrm{Inn} } ( C^{*}( \Z , X, \alpha) ) } { \overline{ \mrm{Inn} }_{0} ( C^{*}( \Z , X, \alpha) ) }$ is a totally disconnect abelian topological group; and

\item The following sequence is exact:
\begin{equation*}
\{ 1 \} \to \frac{ \overline{ \mrm{Inn} } (C^{*}( \Z , X, \alpha)) } { \overline{ \mrm{Inn} }_{0} ( C^{*}( \Z , X, \alpha) ) } \to \frac{ \mrm{Aut} ( C^{*}( \Z , X, \alpha) ) } { \overline{ \mrm{Inn} }_{0} ( C^{*}( \Z , X, \alpha) ) } \to \mrm{ Aut } ( \totalk( C^{*}( \Z , X, \alpha) ) )_{ + , 1} \to \{ 1 \}.
\end{equation*}
\end{enumerate}
\end{theorem}

\begin{proof}
(1) follows from Theorem 4.6 of \cite{linphillipsTR}.  (2) follows from (1) and Corollary 2.5 of \cite{EllRorAuto}.  (3) follows from (1), Theorem \ref{thm:propC}, and Theorem \ref{thm:inn}.  (4) follows from (1) and Theorem \ref{thm:extgrps}.
\end{proof}

We will now show that a large class of \cstar-algebras arising from minimal homeomorphisms of $\mbb{T}^{n}$ satisfies the assumptions of Theorem \ref{thm:minimalDS}.  

\begin{definition}\label{def:furstenberg}
Let $\theta$ be an irrational number, $f_{k}$ be a continuous function from $\mathbb{T}^{k}$ to $\R$, and 
\begin{equation*}
\setof{ d_{ij} \in \Z_{ \geq 0 } }{ i, j}
\end{equation*} 
be a subset of $\Z_{ \geq 0}$.  Define $h_{n}$ from $\mbb{T}^{n}$ to $\mbb{T}^{n}$ to be the inverse of the homeomorphism 
\begin{equation*}
( \zeta_{1} , \zeta_{2} , \dots, \zeta_{n} ) \mapsto ( e^{ 2 \pi i \theta } \zeta_{1} , e^{ 2 \pi i f_{1} ( \zeta_{1} ) } \zeta_{1}^{ d_{12} } \zeta_{2} , e^{2 \pi i f_{2} ( \zeta_{1} , \zeta_{2} ) } \zeta_{1}^{ d_{13} } \zeta_{2}^{ d_{23} } \zeta_{3} , \dots).
\end{equation*}  

The homeomorphism $h_{n}$ will be called a $n$-dimensional Furstenberg transformation.  If $d_{i , i+1}$ is nonzero for all $i$, then by Theorem 2.1 of \cite{furstenberg}, the dynamical system $( \mbb{T}^{n} , h_{n} )$ is minimal.  This implies that the transformation group \cstar-algebra associated to $( \mbb{T}^{n} , h_{n} )$ is simple.  Also, if each $f_{i}$ satisfies a certain Lipschitz property, then Furstenberg in \cite{furstenberg} proved that $( \mbb{T}^{n} , h_{n} )$ has a unique $h_{n}$-invariant probability measure on $\mbb{T}^{n}$.  This implies that the transformation group \cstar-algebra associated to $( \mbb{T}^{n} , h_{n} )$ is simple with a unique tracial state.       
\end{definition}

\begin{proposition}\label{prop:furstenberg}
Consider the dynamical system $( \mbb{T}^{n} , h_{n} )$, where $h_{n}$ is a $n$-dimensional Furstenberg transformation.  Set $A_{n}$ to be the transformation group \cstar-algebra associated to $( \mbb{T}^{n} , h_{n} )$.  Suppose $( \mbb{T}^{n} , h_{n} )$ is minimal with a unique $h_{n}$-invariant probability measure on $\mbb{T}^{n}$.   Then the image of $K_{0} ( A_{n} )$ under the map $\rho_{ A_{n} }$ is dense in $\mrm{Aff} ( T( A_{n} ) )$.
\end{proposition}

\begin{proof}
We follow the computation given in Example 4.9 of \cite{NCPstablerankRSA}.  Define $\alpha$ from $C( \mbb{T}^{n} )$ to $C( \mbb{T}^{n} )$ by $\alpha ( f ) = f \circ h_{n}^{-1}$.  Set
\begin{equation*}
A_{n} = C^{*} ( \Z , \mbb{T}^{n} , h_{n} ) = C^{*} ( \Z , C( \mbb{T}^{n} , \alpha ) ).
\end{equation*}

By the Pimsner-Voiculescu exact sequence \cite{PVseq}, the following sequence 
\begin{equation*}
\xymatrix{
K_{0} ( C ( \mbb{T}^{n} ) ) \ar[r]^{ \id - \alpha_{*}^{-1} } & K_{0} ( C ( \mbb{T}^{n} ) ) \ar[r] & K_{0} ( A_{n} ) \ar[d]^{ \partial } \\
K_{1} ( A_{n} ) \ar[u]^{ \mrm{exp} } 	& K_{1} ( C ( \mbb{T}^{n} ) ) \ar[l] & K_{1} ( C ( \mbb{T}^{n} ) ) \ar[l]^{ \id - \alpha_{*}^{-1} } 
}
\end{equation*}
is exact.  Let $z$ be the unitary in $C( \mbb{T} )$ which sends $\zeta$ to $\zeta$.  Then $K_{1} ( C( \mbb{T}^{n} ) )$ is a finitely generated free abelian group with $[ z \otimes 1 \otimes \cdots \otimes 1 ]$ as one of the generators.  

Note that $\alpha$ is homotopic to the $*$-homomorphism given by $f \mapsto f \circ h^{-1}$ with
\begin{equation*}
h^{-1} ( \zeta_{1} , \zeta_{2} , \dots , \zeta_{n} ) = ( \zeta_{1} , \zeta_{1}^{ d_{12} } \zeta_{2} , \zeta_{1}^{ d_{13} } \zeta_{2}^{ d_{23} } \zeta_{3} , \dots ).
\end{equation*}
Hence, $\alpha ( [ z \otimes 1 \otimes \cdots \otimes 1 ] ) = [ z \otimes 1 \otimes \cdots \otimes 1 ]$.  Therefore, by Proposition 6.1 of \cite{RManzaiflow}, we have that $K_{0} ( A_{n} ) \cong \Z^{m} \oplus G$, where $G$ is a finitely generated torsion group and one of the generators of $K_{0} ( A_{n} )$ is an element $\eta_{0}$ such that $\partial ( \eta_{0} ) = [ z^{-1} \otimes 1 \otimes \cdots \otimes 1]$. 

Let $\tau$ be the unique tracial state of $A_{n}$.  Then $\tau$ is induced by a unique $h_{n}$-invariant probability measure $\mu$ on $\mbb{T}^{n}$.  We compute the image of $\eta_{0}$ under the map $\tau_{*}$ from $K_{0} ( A_{n} )$ to $\R$.  Combining Definition VI.8 and Theorem V.12 and VI.11 of \cite{ExelRotation} to get (notation explained afterwards)
\begin{equation*}
\mrm{exp} ( 2 \pi i \tau_{*} ( \eta_{0} ) ) = R_{ \alpha }^{ \mu } ( [ z^{-1} \otimes 1 \otimes \cdots \otimes 1] ).
\end{equation*}
Here $[ z^{-1} \otimes 1 \otimes \cdots \otimes 1 ]$ now represents the homotopy class of the function 
\begin{equation*}
( \zeta_{1} , \zeta_{2} , \dots, \zeta_{n} ) \mapsto \zeta_{1}^{-1}.
\end{equation*}
Following Definitions VI.3 and VI.5 of \cite{ExelRotation}, $R_{\alpha}^{\mu} ( [ v ] )$ is computed by finding a continuous function $f$ from $\mbb{T}^{n}$ to $\R$ such that 
\begin{equation*}
v ( h_{n}^{-1} ( x ) )^{*} v(x) = e^{ i f(x) }
\end{equation*}
for all $x$ in $\mbb{T}^{n}$.  With $v = z^{-1} \otimes 1 \otimes \cdots \otimes 1$ one can easily check that we may choose the function $f$ to be the constant function $2 \pi \theta$.  Hence, $\mrm{ exp } ( 2 \pi i \tau_{*} ( \eta_{0} ) ) = \exp ( 2 \pi i \theta )$.  Therefore, there exists $k$ in $\Z$ such that $\tau_{*} ( \eta_{0} ) = \theta + k$.  Since $\tau_{*} ( [ 1_{ A_{n} } ] ) = 1$, $\rho_{ A_{n} } ( K_{0} ( A_{n} ) )$ contains $\Z + \theta \Z \subset \R = \mrm{Aff} ( T( A_{n} ) )$.  Therefore, the image of $K_{0} ( A_{n} )$ in $\mrm{Aff} ( T( A_{n} ) )$ under the map $\rho_{ A_{n} }$ is dense in $\mrm{Aff} ( T( A_{n} ) )$.

\end{proof} 

\begin{theorem}\label{thm:furstenberg}
Let $( \mbb{T}^{n} , h_{n} )$ be an $n$-dimensional Furstenberg transformation with a unique $h_{n}$-invariant probability measure on $\mbb{T}^{n}$.  Let $A_{n}$ be the transformation group \cstar-algebra associated to $( \mbb{T}^{n} , h_{n} )$.  Then $\overline{ \mrm{Inn} } ( A_{n} ) = \overline{ \mrm{Inn} }_{0} ( A )$ and $\mrm{Aut} ( A_{n} )$ fits into the following exact sequence:
\begin{equation*}
\xymatrix{
\{ 1 \} \to \overline{ \mrm{Inn} } ( A_{n} )  \to \mrm{Aut} ( A_{n} )  \to \mrm{ Aut }_{ \Lambda } ( \totalk ( A ) )_{ + , 1 } \to \{ 1 \}.
}
\end{equation*}
Consequently, $\overline{ \mrm{Inn} } ( A )$ is a simple topological group.
\end{theorem}

\begin{proof}
By Proposition \ref{prop:furstenberg} and Theorem \ref{thm:minimalDS}, $A_{n}$ is a unital, separable, amenable, simple, tracially AF algebra.  Hence, by Theorem \ref{thm:extgrps}, $\mrm{Aut} ( A_{n} )$ fits into the above short exact sequence.  Note that $K_{0} ( A_{n} )$ is a finitely generated abelian group.  By Theorem \ref{thm:propC} and Theorem \ref{thm:inniso}, $\frac{ \overline{ \mrm{Inn} } ( A_{n} ) } { \overline{ \mrm{Inn}_{0} ( A_{n} ) } }$ is isomorphic to $\frac{ K_{1} ( A_{n} ) }{ H_{ [ 1_{ A_{n} } ] } ( K_{0} ( A_{n} ) , K_{1} ( A_{n} ) ) }$.  By the results of \cite{RManzaiflow}, we have that $K_{0} ( A_{n} )$ is isomorphic to $\Z \oplus G$ in which the isomorphism sends $[ 1_{ A_{n} } ]$ to $(1 , 0 )$.  Hence, $H_{ [ 1_{ A_{n} } ] } ( K_{0} ( A_{n} ) , K_{1} ( A_{n} ) ) = K_{1} ( A_{n} )$.  Therefore, $\overline{ \mrm{Inn} } ( A_{n} ) = \overline{ \mrm{Inn} }_{0} ( A_{n} )$.  Since $A_{n}$ is a unital, simple, tracially AF algebra, $A_{n}$ is a finite \cstar-algebra with real rank zero, has the cancellation property for projections, and $K_{0} ( A_{n} )$ is weakly unperforated.  By Corollary 2.5 of \cite{EllRorAuto}, $\overline{ \mrm{Inn} }_{0} ( A_{n} )$ is a simple topological group.
\end{proof}

In the following example, we will consider a $2$-dimensional Furstenberg transformation.  In particular, we would like to give an explicit description of the group $\mrm{ Aut } ( \totalk(A_{2}) )_{+,1}$.  

\begin{example}
Let $\theta$ be an irrational number in $[ 0 ,1]$, $f_{0}$ from $\mbb{T}$ to $\R$ be a Lipschitz function, and let $d$ be a nonzero integer.  Consider the Furstenberg transformation $h$ from $\mbb{T}^{2}$ to $\mbb{T}^{2}$ defined to be the inverse of the homeomorphism
\begin{equation*}
( \zeta_{1} , \zeta_{2} ) \mapsto ( e^{ 2 \pi i \theta } \zeta_{1}  , e^{ 2 \pi i f_{0} ( \zeta_{1} ) } \zeta_{1}^{d} \zeta_{2} ).
\end{equation*}
If $A$ denotes the transformation group \cstar-algebra associated to $( \mbb{T}^{2} , h )$, then Phillips in \cite{NCPstablerankRSA} (see Example 4.9) showed that $K_{0} ( A )$ is isomorphic to $\Z^{3}$, $K_{1} ( A )$ is isomorphic to $\Z^{3} \oplus \Z / d \Z$, and $K_{0} ( A )_{+}$ can be identified with
\begin{equation*}
\setof{ ( m_{1} , m_{2} , m_{3} ) \in \Z^{3} }{ m_{1} + m_{2} \theta > 0 \ \mrm{or} \ m_{1} = m_{2} = m_{3} = 0}.
\end{equation*}  

Set $G$ to be the following subgroup of $\mrm{GL}_{3} ( \Z )$:
\begin{equation*}
\setof{ 
\left(
\begin{array}{ccc}
1 & 0 & 0 \\
0 & 1 & 0 \\
0 & a  & 1 
\end{array}
\right), \
\left(
\begin{array}{ccc}
1 & 0 & 0 \\
0 & 1 & 0 \\
0 & b  & -1 
\end{array}
\right) }{ a, b \in \Z }.
\end{equation*}
By a simple computation, one can show that $\mrm{Aut} ( K_{0} ( A ) )_{+,1} = G$.  

Since $A$ satisfies the Universal Coefficient Theorem, as noted in \cite{multcoeff} pp. 375, $\mrm{Aut}_{ \Lambda } ( \totalk( A ) )$ can be identified with the groups of units in the ring $\kl ( A , A )$.  Since Dadarlat and Loring's Universal Coefficient Theorem splits,
\begin{equation*}
\kl ( A, A ) \cong \Hom ( K_{*} ( A ) , K_{*} ( A ) ) \oplus \ext ( K_{*} ( A ) , K_{*+1} ( A ) ).
\end{equation*}
The multiplication is induced by the usual action of $\mrm{Hom}$ on $\mrm{Ext}$, which passes to an action on $\mrm{Hom}$ on $\mrm{ext}$, since the pullback and the pushout of a pure extension are pure.  The product of any two elements of $\mrm{ext}$ is zero.  Hence,
\begin{equation*}
\mrm{Aut}_{ \Lambda } ( \totalk ( A ) ) \cong \mrm{ Aut } ( K_{*} ( A ) ) \times \ext ( K_{1} ( A ) , K_{0} ( A ) ).
\end{equation*}
Therefore, we have a bijection 
\begin{align*}
\mrm{Aut}_{ \Lambda } ( \totalk ( A ) )_{+,1} &\cong \mrm{Aut} ( K_{0} ( A ) )_{+,1} \times \mrm{ Aut } ( K_{1} ( A ) ) \times \ext ( K_{1} ( A ) , K_{0} ( A ) ) \\
					&\cong G \times \mrm{Aut} ( \Z^{3} \oplus \Z/ d\Z ) \times \left( \Z / d \Z \right)^{3}.
\end{align*}
The product of two elements is given by 
\begin{equation*}
( \alpha^{0} , \alpha^{1} , x ) \circ ( \beta^{0} , \beta^{1} , y ) = ( \alpha^{0} \beta^{0} , \alpha^{1} \beta^{1} , \alpha^{0} (y) + k_{ \beta^{1} } x )
\end{equation*}  
where $k_{ \beta^{1} }$ is the element of $\Z / d \Z$ such that the projection of $\beta^{1} ( 0 , 1 )$ onto the $\Z / d\Z$ is $k_{ \beta^{1} }$.  Therefore, $\mrm{Aut} ( A )$ fits into the following exact sequence:
\begin{equation*}
\{ 1 \} \to \overline{ \mrm{Inn} } ( A ) \to \mrm{Aut} ( A ) \to G \times \mrm{Aut} ( \Z^{3} \oplus \Z / d \Z ) \times \left( \Z / d \Z \right)^{3} \to \{ 1 \}.
\end{equation*} 
\end{example}

We now discuss the the automorphism group of a simple, unital, AT algebra.  Let $A$ be a unital \cstar-algebra.  Denote the canonical affine map from $T(A)$ to the state space $S ( K_{0} ( A ) )$ of $K_{0} (A)$ by $r_{A}$.  Suppose the canonical map from $U(A) / U_{0} (A)$ to $K_{1} (A)$ is surjective, then there is an embedding 
\begin{equation*}
\ftn{ \lambda_{A} }{ \mrm{Aff} ( T(A) ) / \overline{ \rho_{A} ( K_{0} ( A ) ) } } {  U(A) / CU ( A )  }
\end{equation*}  
which identifies $\mrm{Aff} ( T(A) ) / \overline{ \rho_{A} ( K_{0} ( A ) ) }$ with the connected component of the zero element in 
$U(A) / CU ( A )$ (see \cite{NielsenThomsen}).

For a unital \cstar-algebra $A$, set 
\begin{equation*}
\msf{I}(A) = ( K_{0} ( A ) , K_{0} ( A )_{+} ,  T(A) , U(A) / CU (A)  ).
\end{equation*}
Suppose the canonical map from $U(A) / U_{0} ( A )$ to $K_{1} ( A )$ is surjective.  The group of automorphism of $\msf{I}(A)$ preserving $[1_{A}]$, denoted by $\mrm{Aut} ( \msf{I} ( A ) )_{1}$, consists of a triple $( \alpha_{0} , \alpha_{T}, \beta )$ such that $\alpha_{0}$ is an order automorphism of $( K_{0} (A)  , K_{0} ( A )_{+} )$ with $\alpha_{0} ( [1_{A} ] ) = [ 1_{A} ]$, $\alpha_{T}$ is a affine homeomorphism from $T(B)$ to $T(A)$, and $\beta$ is a contractive group isomorphism from $U ( A ) / CU (A) $ to $U ( B ) /  CU (B) $ which makes the diagram  
\begin{equation*}
\xymatrix{
\mrm{Aff} ( T(A) ) / \overline{ \rho_{A} ( K_{0} ( A ) ) } \ar[r]^(.6){ \lambda_{A} } \ar[d]_{ \widetilde{ \alpha } } & U(A) / CU (A)  \ar[d]^{ \beta } \\
\mrm{Aff} ( T(B) ) / \overline{ \rho_{B} ( K_{0} ( B ) ) } \ar[r]_(.6){ \lambda_{B} } & U(B) / CU ( B )  
}
\end{equation*}
commutative, where $\widetilde{ \alpha }$ is the map induced by $( \alpha_{T} )_{*}$ from $\mrm{Aff} ( T(A) )$ to $\mrm{Aff} ( T(B) )$.

\begin{theorem}\label{thm:simpleAT}
Let $A$ be a unital, simple, AT algebra.  Then $\frac{ \overline{ \mrm{Inn} } ( A ) }{ \overline{ \mrm{Inn} }_{0} ( A ) }$ is a totally disconnected topological group and 
\begin{equation*}
\{ 1 \} \to \frac{ \overline{ \mrm{Inn} } ( A ) }{ \overline{ \mrm{Inn} }_{0} ( A ) } \to \frac{ \overline{ \mrm{Aut} } ( A ) }{ \overline{ \mrm{Inn} }_{0} ( A ) } \to \mrm{ Aut } ( \msf{I}(A) )_{1} \to \{ 1 \}
\end{equation*}
is a short short exact sequence.
\end{theorem}

\begin{proof}
Since $A$ is a unital, simple, AT algebra, by the results of Lin (see Section 7.1 of \cite{Lintracialrank}) we have that $A$ is a tracially AI algebra.  Hence, by Theorem \ref{thm:inn} and Theorem \ref{thm:propC} $\frac{ \overline{ \mrm{Inn} } ( A ) } { \overline{ \mrm{Inn} }_{0} (A) }$ is totally disconnected.    

By Theorem A of \cite{NielsenThomsen}, every element of $\mrm{Aut} ( \msf{I}(A) )_{1}$ lifts to a unital $*$-homomorphism and by Theorem B of \cite{NielsenThomsen}, two automorphisms $\phi$ and $\psi$ of $A$ are approximately unitarily equivalent if and only if $\msf{I}( \phi ) = \msf{I} ( \psi )$.  By these facts and arguing as in Theorem \ref{thm:extgrps}, we get that an element of $\mrm{Aut} ( \msf{I} ( A ) )_{1}$ lifts to an automorphism of $A$.  Hence, the sequence $\mrm{Aut} ( A ) \to \mrm{ Aut } ( \msf{I} ( A ) )_{1} \to \{ 1 \}$ is exact.  By Theorem B of \cite{NielsenThomsen}, we have that the sequence
\begin{equation*}
\{ 1 \} \to \overline{ \mrm{Inn} } ( A ) \to \mrm{ Aut } ( A ) \to \mrm{Aut} ( \msf{I} ( A ) )_{1}
\end{equation*}
is an exact sequence.  It is now clear that the sequence in the theorem is exact.  
\end{proof}

There are lots of examples of \cstar-algebras arising from dynamical systems that are unital, simple, AT algebras.  One example is the irrational rotation algebra.  Other examples are given by Lin and Matui in \cite{linmatui3}.  We described their \cstar-algebras here. Let $( X , \alpha )$ be a dynamical system with $X$ the Cantor set.  For $\lambda$ in $\mbb{T}$, denote the homeomorphism of $\mbb{T}$ which sends $\zeta$ to $\lambda \zeta$ by $R_{ \lambda }$.  Suppose $\xi$ from $X$ to $\mbb{T}$ is a continuous function.   The homeomorphism of $X \times \mbb{T}$ which sends $( x, \zeta )$ to $( \alpha ( x ) , R_{ \xi (x) } ( \zeta ) )$ will be denoted by $\alpha \times R_{ \xi }$.  Lin and Matui in \cite{linmatui3} proved that $C^{*}( \Z , X \times \mbb{T}  , \alpha \times R_{ \xi } )$ is a unital, simple, AT algebra whenever $\alpha \times R_{ \xi }$ is minimal (see Theorem 4.3 of \cite{linmatui3}).  In \cite{linmatui2}, Lin and Matui (Lemma 4.2 of \cite{linmatui2}) showed that $\alpha \times R_{ \xi }$ is minimal if and only if $\alpha$ in minimal and $n [ \xi ] \neq 0$ in $C( X , \mbb{T} ) / \setof{ \eta - \eta \circ \alpha^{-1} }{ \eta \in C( X , \mbb{T} ) }$ for all $n$ in $\N$.  

Note that there exists a canonical factor map from $( X \times \mbb{T} , \alpha \times R_{ \xi } )$ to $( X , \alpha )$.  We say that $\alpha \times R_{ \xi }$ is \emph{rigid} if the canonical factor map from $( X \times \mbb{T} , \alpha \times R_{ \xi } )$ to $( X , \alpha )$ induces an isomorphism between the sets of invariant probability measures.  Theorem 4.3 of \cite{linmatui3} states that if $\alpha \times R_{ \xi }$ is a minimal homeomorphism, then $C^{*} ( \Z , X \times \mbb{T} , \alpha \times R_{ \xi } )$ is a unital, simple, AT algebra with real rank zero if and only if $\alpha \times R_{ \xi }$ is rigid.  They also showed that $\alpha \times R_{ \xi }$ could be minimal but not rigid (see Remark 3.2 of \cite{linmatui2}).        
 
\begin{question}
\begin{enumerate}
\item Is $\overline{ \mrm{Inn} }_{0} ( A )$ a simple topological group for all unital, simple, tracially AI algebra?  
\item Is $\frac{ \overline{ \mrm{Inn} } (A ) } { \overline{ \mrm{Inn} }_{0} (A ) }$ totally disconnected for all simple \cstar-algebras?
\end{enumerate}
\end{question}

\section{acknowledgements}
The first author would like to thank George Elliott and the Fields Institute for their hospitality while he was a post doctoral fellow at the Fields Institute (December 2005 to July 2006), where this research was completed.  The second author would like to thank Fr\'{e}d\'{e}ric Latr\'{e}moli\`{e}re for many stimulating conversations.

\end{document}